\documentclass[reqno,english]{amsart}%
\usepackage{amsfonts,amsmath,latexsym,verbatim,amscd,mathrsfs,color,array}
\usepackage[colorlinks=true]{hyperref}
\usepackage{amsmath,amssymb,amsthm,amsfonts,graphicx,color}
\usepackage{amssymb}
\usepackage{cite}
\usepackage{graphicx}
\usepackage{amsmath}
\usepackage{amsfonts}%
\providecommand{\U}[1]{\protect\rule{.1in}{.1in}}

\numberwithin{equation}{section}

\newtheorem{theorem}{Theorem}[section]
\newtheorem{corollary}{Corollary}[section]
\newtheorem{lemma}{Lemma}[section]
\newtheorem{proposition}{Proposition}[section]
\newtheorem{remark}{Remark}[section]
\newtheorem{definition}{Definition}[section]

\numberwithin{equation}{section}

\newcommand{\bbr}{\mathbb{R}}

\newcommand{\ds}{\displaystyle}
\newcommand{\R}{\mathbb R}

\newcommand{\bd}{\begin{definition}}
	\newcommand{\ed}{\end{definition}}

\newcommand{\br}{\begin{remark}}
	\newcommand{\er}{\end{remark}}

\newcommand{\be}{\begin{equation}}
	\newcommand{\ee}{\end{equation}}

\newcommand{\bc}{\begin{corollary}}
	\newcommand{\ec}{\end{corollary}}

\allowdisplaybreaks[2]
\begin{document}
	\title[Infinitely many positive solutions]{Infinitely many nonradial positive solutions for two-coupled nonlinear Schr\"odinger equations in $\bbr^N$ without symmetry}
	\author[R He]{Rui He}
	\address{\noindent  School of Mathematics, Yunnan Normal University, Kunming, 650500, P.R. China\\
	Yunnan Key Laboratory of Modern Analytical Mathematics and Applications, Kunming, 650500, P.R. China}
	\email{heruiy@ynnu.edu.cn}
	
	\author[X. Liu]{Xiangqing Liu}
\address{\noindent  School of Mathematics, Yunnan Normal University, Kunming, 650500, P.R. China\\
	Yunnan Key Laboratory of Modern Analytical Mathematics and Applications, Kunming, 650500, P.R. China}
	\email{liuxiangqing@ynnu.edu.cn}
	
	\author[Q. Wang]{Qingfang Wang}
	\address{\noindent School of Mathematics and Computer Science, Wuhan Polytechnic University, Wuhan
430079, P. R. China}
	\email{wangqingfang@whpu.edu.cn}
	
	\author[Y. Wu]{Yuanze Wu}
\address{\noindent  School of Mathematics, Yunnan Normal University, Kunming, 650500, P.R. China\\
	Yunnan Key Laboratory of Modern Analytical Mathematics and Applications, Kunming, 650500, P.R. China}
	\email{yuanze.wu@ynnu.edu.cn (Corresponding author.)}
	
	\begin{abstract}
In this paper, we consider the following two-coupled nonlinear
Schr\"odinger equations:
\begin{equation*}
\begin{cases}
&-\Delta u+P(x)u=\mu u^3+\beta uv^2,\quad\text{in }\R^N,\\
&-\Delta v+Q(x)v=\nu v^3+\beta u^2v,\quad\text{in }\R^N,\\
&u,\,v\in H^1(\R^N)
\end{cases}
\end{equation*}
where $N=2,3$, $P(x)$ and $Q(x)$ are potentials,
$\mu$ and $\nu$ are positive coefficients and
$\beta\in\R\setminus\{0\}$ is the coupling parameter.  By developing the intermediate reduction arguments in [del Pino-Wei-Yao, Calc. Var., 2015], we prove the existence of infinitely many nonradial positive solutions to the above system under suitable regularity and decay conditions on $P(x)$ and $Q(x)$, {\bf but without any symmetry assumptions on them}.  The scalar version of our result also gives {\bf a positive answer} to the first part of [Question~3, del Pino-Wei-Yao, Calc. Var., 2015] and our method might be helpful in studying other problems.

\vspace{3mm} \noindent{\bf Keywords:} Coupled Schr\"odinger equations; Infinitely many solutions; Nonradial positive solutions; No symmetric potentials; intermediate reduction arguments.

\vspace{3mm}\noindent {\bf AMS} Subject Classification 2020: 35J50; 35B33; 35B40; 35Q55.%

\end{abstract}
	
	\date{}
	\maketitle

\section{Introduction}
In this paper, we consider the following two-coupled nonlinear
Schr\"odinger equations:
\begin{equation}\label{eq:system}
\begin{cases}
&-\Delta u+P(x)u=\mu u^3+\beta uv^2\quad\text{in }\R^N,\\
&-\Delta v+Q(x)v=\nu v^3+\beta u^2v\quad\text{in }\R^N,
\end{cases}
\end{equation}
where $N=2,3$ is the dimension, $P(x)$ and $Q(x)$ are the trapping potentials,
$\mu$ and $\nu$ are the positive self-interaction coefficients, and
$\beta\in\R\setminus\{0\}$ is the coupling parameter.

\subsection{Background}
Under the standing-wave relations
$\Psi_1(t,x)=e^{-\mathrm i\lambda_1t}u(x)$ and
$\Psi_2(t,x)=e^{-\mathrm i\lambda_2t}v(x)$, where $\mathrm i^2=-1$,
solutions of~\eqref{eq:system} give standing waves for the Gross--Pitaevskii
system
\begin{equation*}\label{eqn9877}
\begin{cases}
\mathrm i\partial_t\Psi_1
=\Delta\Psi_1-(P(x)-\lambda_1)\Psi_1
 +\mu|\Psi_1|^2\Psi_1+\beta|\Psi_2|^2\Psi_1,\\
\mathrm i\partial_t\Psi_2
=\Delta\Psi_2-(Q(x)-\lambda_2)\Psi_2
 +\nu|\Psi_2|^2\Psi_2+\beta|\Psi_1|^2\Psi_2,
\end{cases}
\end{equation*}
with $\Psi_j(t,\cdot)\in H^1(\R^N;\mathbb C)$ for every $t$.
Such systems arise in nonlinear optics and in multicomponent
Bose--Einstein condensates; see
\cite{AkhmedievAnkiewicz1999,EsryEtAl1997,HallEtAl1998,Timmermans1998}.
In optical models, the components describe interacting waves in a nonlinear
medium, while in condensate models, they represent the wave functions of the
different species. The cubic terms describe interactions within and between
components, and the external potentials describe trapping effects. With the
sign convention of~\eqref{eq:system}, $\beta>0$ corresponds to attractive
interaction and $\beta<0$ to repulsive interaction.

According to its important applications in physics, the coupled Schr\"odinger equations have been studied extensively in the literature in the past two decades.  Among which, a particular case is the following equal-frequency case
\begin{equation}\label{eq:system-lim}
\begin{cases}
&-\Delta u+u=\mu u^3+\beta uv^2\quad\text{in }\R^N,\\
&-\Delta v+v=\nu v^3+\beta u^2v\quad\text{in }\R^N,\\
&u(x)\to0\text{ and }v(x)\to0\quad\text{as }|x|\to\infty,
\end{cases}
\end{equation}
since it is easy to see that there exist radial solutions of \eqref{eq:system-lim} as follows
\begin{equation}\label{eqn201}
(w_{1,\beta},w_{2,\beta})
=(c_{1,\beta}w,c_{2,\beta}w),
\end{equation}
provided
\begin{equation}\label{eq:amplitudes}
c_{1,\beta}=\sqrt{\frac{\nu-\beta}{\mu\nu-\beta^2}}\quad\text{and}\quad
c_{2,\beta}=\sqrt{\frac{\mu-\beta}{\mu\nu-\beta^2}}
\end{equation}
for $0<\beta<\min\{\mu,\nu\}$ and $\beta>\max\{\mu,\nu\}$, where $w$ is the positive radial solution of the following equation
\begin{equation}\label{eq:w}
\left\{\aligned
&-\Delta w+w=w^3\quad\text{in }\R^N,\\
&w(0)=\max_{\R^N}w,\,w\in H^1(\R^N),
\endaligned\right.
\end{equation}
whose uniqueness, decay, Morse index one and nondegeneracy are recalled in
\cite[Lemma~3.1]{DelPinoWeiYao2015}.
Moreover, after the work of Lin and Wei~\cite{LinWei2005}, the existence, multiplicity,
and qualitative properties of solutions to autonomous Schr\"odinger systems ($P(x)$ and $Q(x)$ are constants)
have been investigated extensively. We refer to
\cite{AmbrosettiColorado2007,BartschDancerWang2010,BartschWangWei2007,
CLWY26,EFPV2026,LiMiyagakiWang2026,LiuWang08,
MaiaMontefuscoPellacci2006,PWW19,Sirakov2007,SoaveTavares2016,
TV09,WeiYao2012,WeiWu2020} and the references therein.  In particular, very recently,
the weak-coupling part of Sirakov's equal-frequency uniqueness conjecture
is established in \cite{CLWY26}, global bifurcation of nodal
solutions is studied in \cite{LiMiyagakiWang2026} and arbitrarily many
nonradial positive solutions in a strongly competitive limit are constructed
in \cite{EFPV2026}. We also refer to the survey~\cite{T24}
for the connection between elliptic systems and optimal partition problems via a strongly competitive limit.

For non-autonomous systems ($P(x)$ and $Q(x)$ are non-constant functions), it is well-known that trapping potentials affect the location and
shape of positive states, see, for example,
\cite{BML23,DuanGaoGuo2026,IkomaTanaka2011,LinWei2006,
MontefuscoPellacciSquassina2008,PellacciEtAl2024,Pomponio2006} and the references therein.  In particular,
Peng and Wang~\cite{PengWang2013} construct infinitely many nonradial
positive solutions of two different types for \eqref{eq:system} in $\R^3$, under {\bf radial} and
decay assumptions on $P(x)$ and $Q(x)$.
The construction in~\cite{PengWang2013} develops the idea of
Wei and Yan~\cite{WeiYan2010} from the scalar equation to the system.  They use the number of peaks as a large parameter, place the peaks near a regular
polygon and use its radius to balance the potential term against the interaction
between neighboring peaks.  However, unlike the scalar case in \cite{WeiYan2010}, different variational arguments are needed in~\cite{PengWang2013} for the final one-dimensional problem, depending on the coupling and
the sign of the decay coefficients in the leading order. Related reduction arguments are developed in
\cite{ChenPistoia2024,LiWeiWu2024,PistoiaVaira2022,WangYe2025,WangZhao2017},
with different coupling ranges, component numbers, and symmetry assumptions, where we refer \cite{LiWeiWu2024} for the detailed descriptions of this topic.  It is also
worth pointing out that for the scalar equations
\begin{equation}\label{eq:w1}
\left\{\aligned
&-\Delta w+V(x)w=w^3,\quad\text{in }\R^N,\\
&w\in H^1(\R^N),
\endaligned\right.
\end{equation}
a variational approach to the multiplicity of solutions was developed in~\cite{DevillanovaSolimini2012} for radial potential $V(x)$.  However, to our best knowledge, this nice method has not been generalized to the coupled systems.

On the other hand, related multiplicity results for systems like \eqref{eq:system} with non-symmetric potentials
appear in~\cite{AoWangYao2016,KwonLeeLee2022,LiuWangWang2024,
WangXieZhanZhangZhao2015,WangZhou2020}, with their respective structural or
perturbative assumptions.  These results are obtained by applying the methods developed for the scalar multiplicity problem to \eqref{eq:w1} in \cite{AoWei2014,CeramiDevillanovaSolimini2005,
CeramiPassaseoSolimini2013,DevillanovaSolimini2015,
DevillanovaSolimini2020,DelPinoWeiYao2015}, where the perturbative results in~\cite{AoWei2014,CeramiPassaseoSolimini2013,DevillanovaSolimini2015,
DevillanovaSolimini2020}
use a global smallness condition on the potentials.  However, del Pino, Wei and
Yao~\cite{DelPinoWeiYao2015} introduce an intermediate reduction for
algebraically decaying potentials without such a condition.
Roughly speaking, the idea of the intermediate reduction is first
to reduce the final finite-dimensional problem into $\R^2$ by the {\bf reflection symmetry} in higher-dimensional construction, next to solve a large discrete system, which is the discretization of the Jacobi equation of $\mathbb S^1$ by the large number of peaks, by adjusting individual radial and tangential positions with quantitative bounds, and finally to kill the multiplier in the kernel of this large discrete system, which comes from the rotational invariance of the reference polygon, by variational arguments for the action of $SO(2)$.  The application of this intermediate reduction to
a cubic two-component system in $\R^2$ is achieved in \cite{AoWangYao2016} with a {\bf smallness} condition on the coupling parameter so that the system is treated as a perturbation of two independent scalar equations.

\subsection{Main results}
The purpose of this paper is to remove the reflection symmetry on potentials and smallness condition on the coupling parameter in constructing infinitely many nonradial positive solutions of systems. For this construction, we impose
the following assumptions:
\begin{enumerate}
\item[$(C_1)$] $P(x),Q(x)\in L^\infty(\R^N)\cap C^{0,\alpha}_{\mathrm{loc}}(\R^N)$
for some $\alpha\in(0,1)$,
\begin{eqnarray*}
P_0:=\inf_{x\in\R^N}P(x)>0\quad\text{and}\quad Q_0:=\inf_{x\in\R^N}Q(x)>0.
\end{eqnarray*}
\item[$(C_2)$] Uniformly as $|x|\to\infty$,
\begin{equation*}\label{eq:tails}
P(x)=1+\frac{a}{|x|^{m}}+O\left(\frac{1}{|x|^{m+\sigma}}\right)\quad\text{and}\quad
Q(x)=1+\frac{b}{|x|^{n}}+O\left(\frac{1}{|x|^{n+\sigma}}\right),
\end{equation*}
where $a,b\in\R$, $m_0=\min\{m,n\}>2$ and $\sigma>2$.
\end{enumerate}
Without loss of generality, we assume
$0<\mu\le\nu$.    Then
for $\beta\in(0,\mu)\cup(\nu,\infty)$,
by \cite[Theorem~1]{Sirakov2007} and \cite[Theorem~1.1]{CLWY26}, the only positive solution of system \eqref{eq:system-lim} is given by \eqref{eqn201}.  Its nondegeneracy follows from the weighted spectral calculation associated
with this profile, see \cite[Lemma~3.1]{BartschDancerWang2010} and
Proposition~\ref{mn:profilekernel} below.
If $\mu<\nu$ and $\mu\le\beta\le\nu$, then by \cite[Theorem~1]{Sirakov2007}, no positive finite-energy
solution of the equal-frequency autonomous system~\eqref{eq:system-lim} exists.  At $\mu=\nu=\beta$, instead, by \cite[Remark~2]{Sirakov2007}, the positive family is
\begin{equation*}\label{eqn2011}
(w_{1,\theta},w_{2,\theta})
=\frac1{\sqrt{\mu}}(w\cos\theta,\,w\sin\theta\,)
\end{equation*}
for all $0<\theta<\frac\pi2$.
We call this the resonant case considered in the present paper. The term
nonresonant below refers to the two positive ranges in~\eqref{eq:amplitudes}.

\vskip0.12in

Let $K\in\mathbb{N}$.
For $j=0,\ldots,K-1$ we denote $\theta_j=2\pi j/K$ and introduce
\begin{equation*}\label{eq:frame}
\pmb n_j=\left\{\aligned
&(\cos\theta_j,\sin\theta_j),\quad N=2,\\
&(\cos\theta_j,\sin\theta_j,0),\quad N=3
\endaligned\right.
\quad\text{and}\quad
\pmb t_j=
\left\{\aligned
&(-\sin\theta_j,\cos\theta_j),\quad N=2,\\
&(-\sin\theta_j,\cos\theta_j,0),\quad N=3.
\endaligned\right.
\end{equation*}
We also define
\begin{equation}\label{eq:synccenters}
Q_j=\begin{cases}
\mathcal A((R+f_j)\pmb n_j+g_j\pmb t_j),&N=2,\\
\mathcal A((R+f_j)\pmb n_j+g_j\pmb t_j+z_j\pmb e_3),&N=3,
\end{cases}
\end{equation}
where $R>0$ is chosen below and $\mathcal A\in SO(N)$. We write the configuration as
\begin{equation}\label{eq:configuration}
\overline{\mathbf q}=\left\{\aligned
&(\mathbf f,\mathbf g),\quad N=2,\\
&(\mathbf f,\mathbf g,\mathbf z),\quad N=3,
\endaligned\right.
\end{equation}
where $\mathbf f=(f_0, f_1, \cdots, f_{K-1})$, $\mathbf g=(g_0, g_1, \cdots, g_{K-1})$ and $\mathbf z=(z_0, z_1, \cdots, z_{K-1})$.  Moreover, we periodically extend $\mathbf f$, $\mathbf g$ and $\mathbf z$ to $\mathbb{Z}$ to be cyclic sequences.
For a cyclic sequence $\mathbf a=(a_0, a_1, \cdots, a_{K-1})$ and $h=2\pi/K$, we introduce its discrete norm
\begin{equation}\label{eq:discretenorm}
\|\mathbf a\|_*=
\|\mathbf a\|_\infty+\|\dot{\mathbf a}\|_\infty+\|\ddot{\mathbf a}\|_\infty,
\end{equation}
as in \cite{DelPinoWeiYao2015}, where
\begin{equation*}\label{eq:difference}
\dot a_j=\frac{a_{j+1}-a_j}{h}\quad\text{and}\quad
\ddot a_j=\frac{a_{j+1}-2a_j+a_{j-1}}{h^2}
\end{equation*}
are its forward first difference and centered second difference, respectively.
The same norm of the configuration $\overline{\mathbf q}$ is defined componentwise.  The associated energy of \eqref{eq:system} is given by
\begin{equation*}
\begin{aligned}
I_\beta(u,v)={}&\frac12\int_{\R^N}
 (|\nabla u|^2+P(x)u^2+|\nabla v|^2+Q(x)v^2)\, dx\\
&-\frac14\int_{\R^N}(\mu u^4+\nu v^4+2\beta u^2v^2)\, dx.
\end{aligned}
\end{equation*}
Now, our result for the nonresonant case can be stated as follows.
\begin{theorem}\label{Nonresonant}
Let $N=2,3$, $0<\mu\le\nu$ and
$\beta\in(0,\mu)\cup(\nu,\infty)$ be fixed. Assume $(C_1)$--$(C_2)$ and
one of the following sign conditions hold:
\begin{equation}\label{eq:signcondition}
\begin{cases}
a>0,&m<n,\\
b>0,&n<m,\\
a c_{1,\beta}^2+b c_{2,\beta}^2>0,&m=n.
\end{cases}
\end{equation}
Then, for every sufficiently large $K\in\mathbb{N}$, there exists $R_K$ and $\overline{\mathbf q}$ such that system
\eqref{eq:system} has a positive solution
\begin{equation*}\label{eq:syncnrprofile}
\binom{u_K}{v_K}
=\sum_{j=0}^{K-1}
\binom{w_{1,\beta}(x-Q_j)}{w_{2,\beta}(x-Q_j)}
+\binom{\phi_K}{\psi_K},
\end{equation*}
where $\|(\phi_K,\psi_K)\|_\infty+\|(\phi_K,\psi_K)\|_{H^1(\R^N)}\rightarrow0,\,\,\hbox{as }\,K\rightarrow+\infty$.
Moreover,
\begin{equation*}\label{eq:syncenergy}
I_\beta(u_K,v_K)
=K(c_{1,\beta}^2+c_{2,\beta}^2)I(w)+o(1).
\end{equation*}
Thus, \eqref{eq:system} has infinitely many nonradial positive solutions with unbounded energy.
\end{theorem}

\begin{remark}\label{rmk01}
The radius $R_K$ satisfies
\begin{eqnarray}\label{eq:R-asymptotic}
R_K=\frac{m_0}{2\pi}K\log K
 +\frac{m_0-(N-3)/2}{2\pi}K\log(m_0\log K)+O(K).
\end{eqnarray}
\end{remark}
\begin{remark}\label{rmk02}
The configuration $\overline{\mathbf q}$ satisfies the following orthogonal conditions
\begin{equation}\label{eq:spatialconditions}
\begin{cases}
\displaystyle\sum_{j=0}^{K-1}g_j=0,&N=2,\\[1mm]
\displaystyle\sum_{j=0}^{K-1}g_j=\sum_{j=0}^{K-1}z_j\cos\theta_j=\sum_{j=0}^{K-1}z_j\sin\theta_j=0,&N=3.
\end{cases}
\end{equation}
Moreover, we also have the following controls
\begin{equation*}
\|\overline{\mathbf q}\|_*\lesssim\delta_K,\quad
\|(\phi_K,\psi_K)\|_\infty\lesssim R_K^{-m_0}\quad\text{and}\quad
\|(\phi_K,\psi_K)\|_{H^1(\R^N)}
\lesssim \sqrt K\,R_K^{-m_0},
\end{equation*}
with
\begin{equation}\label{eq:deltaK}
\delta_K=R_K^{2-\min\{\sigma,m_0\}}+R_K^{-1}+R_K^{2-(1-\eta)m_0}
\end{equation}
for $\eta>0$ sufficiently small that $(1-\eta)m_0>2$.
\end{remark}
As pointed out above,  the limiting system~\eqref{eq:system-lim} has an additional component-rotation
symmetry at resonance. Thus, to achieve a similar construction,  a joint symmetry of the
potentials $P(x)$ and $Q(x)$ is imposed.
\begin{theorem}\label{Resonant}
Let $N=2,3$ and $\mu=\nu=\beta>0$. Assume $(C_1)-(C_2)$ with $m=n>2,\,a=b>0$ and
 $\sigma>2$.
Assume in addition that $P(x)=Q(-x)$ for all $x\in\R^N$.
Then, for every sufficiently large even integer $K$, there exists $R_K$ and $\overline{\mathbf q}$ such that system~\eqref{eq:system}
has a positive solution
\begin{equation*}\label{eq:syncresprofile}
\binom{u_K}{v_K}
=\frac1{\sqrt{\mu}}\sum\limits_{j=0}^{K-1}
\binom{w(x-Q_j)\cos(\pi/4+a_j)}{w(x-Q_j)\sin(\pi/4+a_j)}
+\binom{\phi_K}{\psi_K},
\end{equation*}
satisfying $v_K(x)=u_K(-x)$.  Moreover, the energy satisfies $I_\mu(u_K,v_K)=\frac{KI(w)}{\mu}+o(1)$.  Thus, \eqref{eq:system} has infinitely many nonradial positive solutions with unbounded energy.
\end{theorem}

\begin{remark}
The radius $R_K$ and the configuration $\overline{\mathbf q}$ satisfies the same conditions in Remarks~\ref{rmk01} and \ref{rmk02}, respectively.  Moreover, the remainders satisfy
\[
 \|(\phi_K,\psi_K)\|_\infty\lesssim R_K^{-m}\quad\text{and}\quad
 \|(\phi_K,\psi_K)\|_{H^1}
       \lesssim \sqrt K R_K^{-m},
\]
while the angles $\mathbf a$ satisfy
$\sum_{j=0}^{K-1}a_j=0$ and $\|\mathbf a\|_*\lesssim\delta_K/d_K$,
where $d_K=2R_K\sin(\pi/K)$ and $m_0=m$ in $\delta_K$ given in Theorem~\ref{Nonresonant}.
\end{remark}
\begin{remark}
Clearly, if $P\not\equiv Q$, then $u_K\not\equiv v_K$.
\end{remark}
\subsection{Sketch of the proof}
The proofs of Theorems~\ref{Nonresonant} and~\ref{Resonant} develop the
intermediate reduction of~\cite{AoWangYao2016,DelPinoWeiYao2015} in dimension three.  More precisely, we use the full space in dimension three by adding
$z_j\pmb e_3$ at every center $Q_j$ and replacing the action of $SO(2)$ by $SO(3)$ on the polygon.
\vskip0.12in

\subsubsection{The nonresonant case}\quad
At $(w_{1,\beta},w_{2,\beta})$, the linearization of \eqref{eq:system-lim} is
\begin{equation}\label{eq:general-linearization}
\mathcal L_\beta=(-\Delta+1)I-w^2
\begin{pmatrix}
3\mu c_{1,\beta}^2+\beta c_{2,\beta}^2&2\beta c_{1,\beta}c_{2,\beta}\\
2\beta c_{1,\beta}c_{2,\beta}&3\nu c_{2,\beta}^2+\beta c_{1,\beta}^2
\end{pmatrix},
\end{equation}
where the matrix
\begin{eqnarray}\label{Idmatrix}
I=\begin{pmatrix}
1&0\\
0&1
\end{pmatrix}.
\end{eqnarray}
Its kernel in the nonresonant case is given by
\begin{equation*}\label{eq:common-translation-kernel}
\ker\mathcal L_\beta
=\operatorname{span}\left\{
\binom{\partial_lw_{1,\beta}}{\partial_lw_{2,\beta}}\mid 1\le l\le N
\right\},
\end{equation*}
which can be found in \cite[Lemma~3.1]{BartschDancerWang2010}.  We recall it in Proposition~\ref{mn:profilekernel}.
For $K$ peaks $(w_{1,\beta}(x-Q_j),w_{2,\beta}(x-Q_j))$, these single-profile kernels yield $NK$ local approximate
kernel directions for the linearization of \eqref{eq:system} at the ansatz
\begin{eqnarray*}
(U,V)=\sum_{j=0}^{K-1}(w_{1,\beta}(x-Q_j),w_{2,\beta}(x-Q_j))
\end{eqnarray*}
in the nonresonant case.  After solving the projected equation by the standard fixed-point argument, the
remaining task is to use the variational argument to make its multipliers in the approximate
kernels vanish, as in \cite{AoWangYao2016,ChenPistoia2024,LiWeiWu2024,PistoiaVaira2022,WangYe2025,WangZhao2017}.  To achieve this goal, it is well known that the principal part of
the common translation projection is obtained from the following identity
\begin{eqnarray}\label{RedPro}
&&\int(\mathcal E_1\nabla w_{1,\beta}(\cdot-Q_k)
 +\mathcal E_2\nabla w_{2,\beta}(\cdot-Q_k))\, dx\notag\\
&=&\int(c_{1,\beta}^2(P(x)-1)+c_{2,\beta}^2(Q(x)-1))W\nabla w_k\, dx\notag\\
&&-(c_{1,\beta}^2+c_{2,\beta}^2)
  \int(W^3-\textstyle\sum_{j=0}^{K-1}w_j^3)\nabla w_k\, dx,
\end{eqnarray}
since we have good control on the remainders in any reasonable sense by adapting the system estimates in \cite{LiWeiWu2024} without any smallness assumption on the coupling parameter,
where
\begin{eqnarray*}
\left\{\aligned
&\mathcal E_1=c_{1,\beta}\bigl((P(x)-1)W-(W^3-\textstyle\sum_{j=0}^{K-1}w_j^3)\bigr),\\
&\mathcal E_2=c_{2,\beta}\bigl((Q(x)-1)W-(W^3-\textstyle\sum_{j=0}^{K-1}w_j^3)\bigr)
\endaligned\right.
\end{eqnarray*}
with $w_j=w(\cdot-Q_j)$ and $W=\sum_{j=0}^{K-1}w_j$.

\vskip0.12in

More precisely, estimating the projected correction and the non-neighbor
interactions gives the reduced force equation
\begin{eqnarray}\label{s:eq:neighborforce}
0&=&B(|Q_k|)\frac{Q_k}{|Q_k|}
+\sum_{j=k\pm1}\Psi(|Q_j-Q_k|)\frac{Q_j-Q_k}{|Q_j-Q_k|}+r_k\notag\\
&:=&F_k+r_k
\end{eqnarray}
for $0\leq k\leq K-1$, where $\Psi(s)=-(c_{1,\beta}^2+c_{2,\beta}^2)H_w'(s)$ with
\begin{equation}\label{eq:interactions-general}
H_w(s)=\int_{\R^N}w^3(x)w(x-s\pmb e_1)\,dx,\,\,\,\,\,
H_w'(s)=-\int_{\R^N}w^3(x)\partial_1w(x-s\pmb e_1)\,dx
\end{equation}
and
\begin{equation}\label{eq:potential-coefficients}
B(r)=a_0m r^{-m-1}+b_0n r^{-n-1}
\end{equation}
with $a_0=\frac a2\ds\int_{\R^N}w_{1,\beta}^2\, dx$ and $b_0=\frac b2\ds\int_{\R^N}w_{2,\beta}^2\, dx$.  As in \cite{AoWangYao2016,ChenPistoia2024,LiWeiWu2024,PistoiaVaira2022,WangYe2025,WangZhao2017}, the balance condition is therefore given by
\begin{equation}\label{eq:balance}
B(R_K)=2\sin\frac\pi K\,\Psi(d_K)\quad\text{with}\quad d_K=2R_K\sin\frac\pi K.
\end{equation}
Indeed, by the system estimates in \cite{LiWeiWu2024} again, we have
$\max_{0\leq k \leq K-1}|r_k|=o(R_K^{-m_0-2})$ in \eqref{s:eq:neighborforce}, not just for \eqref{RedPro}.  Thus, the principal force of the reduced problem is ${\bf F}=(F_0, F_1, \cdots, F_{K-1})$.

\vskip0.12in

Let $s=\sin(\pi/K)$, $c=\cos(\pi/K)$ and
\begin{equation}\label{eq:D-nueff}
D_K=-\frac{d_K\Psi'(d_K)}{\Psi(d_K)},\quad
\nu_{*,K}=-\frac{R_KB'(R_K)}{B(R_K)}\quad\text{and}\quad b_{*,K}:=\frac{B(R_K)}{R_K}.
\end{equation}
The linearization of the principal part $F_k$ at the polygon, divided by $b_{*,K}$ and
written in the local frames 
\begin{eqnarray*}
\left\{\aligned
&(\pmb n_j, \pmb t_j),\quad N=2,\\
&(\pmb n_j, \pmb t_j, \pmb e_3),\quad N=3,
\endaligned\right.
\end{eqnarray*}
is $T_N\overline{\mathbf q}$, where the matrix
\begin{equation}\label{s:eq:T3}
T_N=\left\{\aligned
&T_\parallel,\quad N=2,\\
&\begin{pmatrix}T_\parallel&0\\0&T_\perp\end{pmatrix},\quad N=3,
\endaligned\right.
\end{equation}
for
\begin{equation}\label{eq:matrixblocks}
T_\parallel=\begin{pmatrix}
a_{e,K}A_2+(D_K-\nu_{*,K})I&b_{e,K}A_1\\
-b_{e,K}A_1&c_{e,K}A_2
\end{pmatrix}\quad\text{and}\quad
T_\perp=I+\frac{A_2}{4s^2},
\end{equation}
with $I$ the identity matrix given by \eqref{Idmatrix},
\begin{equation}\label{eq:matrixcoefficients}
a_{e,K}=\frac{c^2+D_Ks^2}{4s^2},\quad
b_{e,K}=\frac{(D_K-1)c}{4s}\quad\text{and}\quad
c_{e,K}=-\frac{D_Kc^2+s^2}{4s^2},
\end{equation}
and $A_1$ and $A_2$ are the circulant matrices
\begin{equation}\label{eq:A1A2-matrices}
\begin{aligned}
A_1=\operatorname{circ}(0,1,0,\ldots,0,-1)\quad\text{and}\quad
A_2=\operatorname{circ}(-2,1,0,\ldots,0,1),
\end{aligned}
\end{equation}
whose actions on a cyclic sequence are given by
\begin{equation*}\label{eq:A1A2}
\begin{aligned}
(A_1a)_j=a_{j+1}-a_{j-1}\quad\text{and}\quad(A_2a)_j=a_{j+1}-2a_j+a_{j-1}.
\end{aligned}
\end{equation*}
The matrix $T_\parallel$ is the one appearing in
\cite[Section~5 and Appendix~A]{DelPinoWeiYao2015}. However, to remove the reflection symmetry for $x_3$ in the original intermediate reduction in \cite{DelPinoWeiYao2015}, we introduce the perturbation $z_j\pmb e_3$ at every center $Q_j$, which generates the matrix $T_{\perp}$ for $N=3$.

For sufficiently large $K$, it is shown in \cite[Section~5 and Appendix~A]{DelPinoWeiYao2015} that
$\ker T_\parallel=\operatorname{span}\{(\mathbf0,\mathbf1)\}$.
It is worth pointing out that this vector is just the tangent to the $SO(2)$-orbit of the polygon
in local coordinates, divided by $R$. In dimension three, we prove that $\ker T_3=\operatorname{span}\{\mathbf k_1,\mathbf k_2,\mathbf k_3\}$, where
\begin{equation}\label{eq:rotationkernels}
\begin{gathered}
\mathbf k_1=(\mathbf0,\mathbf1,\mathbf0),\qquad
\mathbf k_2=(\mathbf0,\mathbf0,(\sin\theta_j)_j),\qquad
\mathbf k_3=(\mathbf0,\mathbf0,(-\cos\theta_j)_j).
\end{gathered}
\end{equation}
Moreover, we observe that these three vectors span the tangent space to the rotation orbit of the polygon after division by $R$.

\vskip0.12in

With a clear understanding of the kernel of $T_N$, we can carry out a secondary reduction on \eqref{s:eq:neighborforce} under the orthogonal conditions \eqref{eq:spatialconditions} on the configuration $\overline{\mathbf q}$, as in \cite{DelPinoWeiYao2015}.  The
remaining equations therefore are the projections of the principal force ${\bf F}$ on the kernel of $T_N$ given by
\begin{equation*}\label{eq:globalmultipliers}
\gamma_a(\mathcal A)=0\quad\text{for}\quad
1\le a\le\frac{N(N-1)}2,
\end{equation*}
which is equivalent to $\nabla_{\mathcal A}F=0$, where 
$$
F(\mathcal{A})=I_\beta\left(u_K,v_K\right)=I_\beta\left(\sum\limits_{j=0}^{K-1}w_{1,\beta}(x-Q_j)+\phi_K,\sum\limits_{j=0}^{K-1}w_{2,\beta}(x-Q_j)+\psi_K\right).
$$
Compactness of $SO(N)$ supplies such a critical point, see for example, \cite{Lee2013}.

\vskip0.12in

\subsubsection{The resonant case}\quad
In the resonant case, with the notations
\[
e(\vartheta)=\binom{\cos\vartheta}{\sin\vartheta},\quad
J=\begin{pmatrix}0&-1\\1&0\end{pmatrix}\quad\text{and}\quad
\vartheta_j=\vartheta+a_j,
\]
we choose the following ansatz
\begin{align*}
(U,V)=\left(\frac{1}{\sqrt{\mu}}\sum_{j=0}^{K-1}w_j\cos(\vartheta+a_j),\frac{1}{\sqrt{\mu}}\sum_{j=0}^{K-1}w_j\sin(\vartheta+a_j)\right).
\end{align*}
Its local approximate kernel consists of
\begin{equation*}\label{eq:reskernel}
\frac1{\sqrt{\mu}}e(\vartheta_j)\partial_lw_j
\text{ for }1\le l\le N\quad\text{and}\quad
Y_j=\frac1{\sqrt{\mu}}w_jJe(\vartheta_j).
\end{equation*}

Thus the first projected problem includes $(N+1)K$ local directions in the kernel.
The computation on the corrected energy yields the coupled leading position and angle system
\begin{equation}\label{r:eq:angleforce}
\left\{\begin{aligned}
0={}&A_0m|Q_j|^{-m-2}Q_j
 +\frac1{\mu}\sum_{i=j\pm1}
 \Psi_w(|Q_i-Q_j|)\frac{Q_i-Q_j}{|Q_i-Q_j|}\cos(a_i-a_j)
 +r_{p,j},\\
0={}&\frac{q_j^*-p_j^*}{2\mu}\sin(2\vartheta_j)
 +\frac1{\mu}\sum_{i=j\pm1}H_w(|Q_i-Q_j|)\sin(a_j-a_i)
 +r_{{\bf a}, j},
\end{aligned}\right.
\end{equation}
where $\Psi_w(s)=-H_w'(s)$, $A_0=\frac{a}{2\mu}\ds\int_{\R^N}w^2\, dx$,
\begin{eqnarray*}
p_j^*=\int_{\R^N}(P(x)-1)w_j^2\, dx\quad\text{and}\quad
q_j^*=\int_{\R^N}(Q(x)-1)w_j^2\, dx.
\end{eqnarray*}
The two remainders $r_{p,j}$ and $r_{{\bf a}, j}$ contain all corrections and the non-neighbor terms.  Moreover,
the local-angle equations are nonlinearly coupled with the position equations.

\vskip0.12in

A crucial point in solving \eqref{r:eq:angleforce} is that its linearization decouples into
\begin{eqnarray*}
\text{diag}\left(T_N, -\frac{H_w(d)}\mu A_2+\frac{q_j^*-p_j^*}{\mu}\cos2\vartheta\right).
\end{eqnarray*}
Thus, we can first solve the
position range equations under the orthogonal condition \eqref{eq:spatialconditions} and the
relative-angle range equations under the orthogonal condition $\sum_{j=0}^{K-1}a_j=0$.
The remaining parameters are $\mathcal A\in SO(N)$ and the common angle
$\vartheta$ in a fixed compact subinterval of $(0,\pi/2)$.
Since the resulting functions $\overline{\mathbf q}(\mathcal A,\vartheta)$
and $\mathbf a(\mathcal A,\vartheta)$ are $C^1$ by the implicit function theorem, we can use the variational argument to kill the multipliers in the kernel.  The remaining energy is
\[
 F(\mathcal A,\vartheta)=
 \mathcal J(\mathcal A,\overline{\mathbf q}(\mathcal A,\vartheta),
                  \vartheta,\mathbf a(\mathcal A,\vartheta))\quad\text{with}\quad \mathcal J=I_\mu(u_K,v_K).
\]
Under the additional assumption $Q(x)=P(-x)$, the transformation
$\mathscr S(u,v)(x)=(v(-x),u(-x))$ preserves the system and its energy.  Thus, the half-period relabeling of an even polygon (see Step.~1 in the proof of Theorem~\ref{Resonant} for details), together with uniqueness in the first projected problem and in the joint range equations, implies
\begin{equation*}\label{eq:exact-angle-symmetry}
F(\mathcal A,\vartheta)=F(\mathcal A,\pi/2-\vartheta)\quad\text{and}\quad
\partial_\vartheta F(\mathcal A,\pi/4)=0.
\end{equation*}
This finally solves the common-angle equation.  Choosing a critical point of $\mathcal A\mapsto F(\mathcal A,\pi/4)$ on
$SO(N)$ then kills the multipliers in the full reduced system, and the uniqueness
relations also yield $v_K(x)=u_K(-x)$.

\subsection{Further remarks}
We can also use the same spirit to construct infinitely many nonradial positive solutions in higher dimensions.  To clearly explain the construction, we move to the following scalar equation
\begin{equation}\label{eq:scalar}
\begin{cases}
&-\Delta u+V(x)u=\mu u^p\quad\text{in }\R^N,\\
&u(x)\to0\quad\text{as }|x|\to\infty,
\end{cases}
\end{equation}
where $N\geq3$, $1<p<\frac{N+2}{N-2}$ and the potential $V(x)$ satisfies
\begin{enumerate}
\item[$(D_1)$]$V(x)\in L^\infty(\R^N)\cap C^{0,\alpha}_{\mathrm{loc}}(\R^N)$
for some $\alpha\in(0,1)$ and
$V_0:=\inf\limits_{x\in\R^N}V(x)>0$.
\item[$(D_2)$] Uniformly as $|x|\to\infty$,
$V(x)=1+\frac{a}{|x|^{m_0}}+O\left(\frac{1}{|x|^{m_0+\sigma}}\right)$,
where $a>0$, $m_0\min\left\{1,\frac{p-1}{2}\right\}>2$ and $\sigma>2$.
\end{enumerate}
We define
\begin{equation}\label{eq:scalarcenters}
Q_{j,s}=
\mathcal A\left((R+f_j)\pmb n_j+g_j\pmb t_j+\sum_{l=3}^{N}z_{j,l}\pmb e_l\right)
\end{equation}
with $\overline{\mathbf q}_{s}=(\mathbf f,\mathbf g,\mathbf z_3, \cdots, \mathbf z_N)$.  The energy of \eqref{eq:scalar} is given by
\begin{eqnarray*}
\begin{aligned}
H(u)={}&\frac12\int_{\R^N}
 (|\nabla u|^2+V(x)u^2)\, dx-\frac{\mu}{p+1}\int_{\R^N}|u|^{p+1}\, dx.
\end{aligned}
\end{eqnarray*}
Now our results for \eqref{eq:scalar} can be stated as follows.
\begin{theorem}\label{Scalar}
Let $N\geq3$, $\mu>0$ and $(D_1)$--$(D_2)$ hold.
Then, for every sufficiently large $K\in\mathbb{N}$, there exists a radius $R_K$ satisfying \eqref{eq:R-asymptotic}
such that equation
\eqref{eq:scalar} has a positive solution
\begin{equation*}\label{eq:syncnrprofile-scalar}
u_K
=\sum_{j=0}^{K-1}
w(x-Q_{j,s})+\phi_K,
\end{equation*}
where the configuration $\overline{\mathbf q}_{s}$ satisfies the orthogonal conditions
\begin{equation}\label{eq:spatialconditions-scalar}
\displaystyle\sum_{j=0}^{K-1}g_j=\sum_{j=0}^{K-1}z_{j,l}\cos\theta_j=\sum_{j=0}^{K-1}z_{j,l}\sin\theta_j=0
\end{equation}
for all $3\leq l\leq N$ with the estimates
\begin{equation*}\label{eq:syncnrestimates}
\|\overline{\mathbf q}_{s}\|_*\lesssim\delta_K,\quad
\|\phi_K\|_\infty\lesssim R_K^{-m_0},\quad
\|\phi_K\|_{H^1(\R^N)}\lesssim \sqrt K\,R_K^{-m_0}
\end{equation*}
and
\begin{equation*}
\delta_K=R_K^{2-\min\left\{\sigma,\min\left\{1,\frac{p-1}{2}\right\}m_0\right\}}+R_K^{-1}+R_K^{2-(1-\eta)m_0}
\end{equation*}
for $\eta>0$ sufficiently small that $(1-\eta)m_0>2$.  Moreover,
$H(u_K)=KI(w)+o(1)$.
Thus, \eqref{eq:scalar} has infinitely many nonradial positive solutions with unbounded energy.
\end{theorem}

\begin{remark}
Theorem~\ref{Scalar} is the scalar version of the main theorems~\ref{Nonresonant} and \ref{Resonant}, which also gives a positive answer to the first part of
\cite[Question~3]{DelPinoWeiYao2015}, since we remove the reflection symmetry on the potential $V(x)=V(x_1,x_2,x')$ for $x'\in\mathbb{R}^{N-2}$.
\end{remark}

The proof of Theorem~\ref{Scalar} is achieved by following the arguments and calculations in \cite{DelPinoWeiYao2015} with two major modifications.  The first modification is to use the full approximate kernel
$\frac1{\sqrt{\mu}}\partial_lw_j
\text{ for }1\le l\le N$
in the first projected problem, according to the new perturbations $\sum_{l=3}^{N}z_{j,l}\pmb e_l$.  The second modification is to replace $T_\parallel$ by
\begin{eqnarray*}
\text{diag}(T_\parallel, \underbrace{T_{\perp}, \cdots, T_{\perp}}_{N-2})
\end{eqnarray*}
to serve as the linearization of the principal force to the reduced problem.  Similar to \eqref{eq:rotationkernels}, its kernel is given by
\begin{equation*}\label{eq:rotationkernels-N}
\left\{\aligned
&\mathbf k_1=(\mathbf0,\mathbf1,\underbrace{\mathbf0,\cdots,\mathbf0}_{N-2}),\\
&\mathbf k_{1,l}=(\underbrace{\mathbf0,\cdots,\mathbf0}_{l-1},(\sin\theta_j)_j,\underbrace{\mathbf0,\cdots,\mathbf0}_{N-l}),\\
&\mathbf k_{2,l}=(\underbrace{\mathbf0,\cdots,\mathbf0}_{l-1},(-\cos\theta_j)_j,\underbrace{\mathbf0,\cdots,\mathbf0}_{N-l})
\endaligned\right.
\end{equation*}
for $3\leq l\leq N$.
Clearly, $\mathbf k_1$ is the tangent vector of the rotation orbit
of the polygon in $\text{span}\{\pmb e_1, \pmb e_2\}$ after division by $R$, while $\mathbf k_{1,l}$ and $\mathbf k_{2,l}$ are the tangent vectors of the rotation orbit of
the polygon in $\text{span}\{\pmb e_1, \pmb e_l\}$ and $\text{span}\{\pmb e_2, \pmb e_l\}$, respectively, after division by $R$.  Moreover, it is also well-known (\cite{MilnorStasheff1974}) that the tangent space of $SO(N)/SO(N-2)$ at the identity is given by
\begin{eqnarray*}
\text{span}_{3\leq l\leq N}\{\mathbf k_1, \mathbf k_{1,l}, \mathbf k_{2,l}\}.
\end{eqnarray*}
Thus, by restricting $\mathcal A\in SO(N)/SO(N-2)$ in \eqref{eq:scalarcenters}, we can kill all multipliers in the projected problem by variational arguments as for \eqref{eq:system} in the nonresonant case.

\subsection{Organization and notation}
Section~\ref{sec:pre-munu} collects the profiles, spectra, interaction
estimates and projected equations. Section~\ref{sec:nonres-munu} proves
Theorem~\ref{Nonresonant}, with the exact cyclic matrix, its inverse,
the nonlinear position equations, and the final rotational argument.
Section~\ref{sec:res-munu} proves Theorem~\ref{Resonant}, by solving the
position and relative-angle equations jointly and then using the
inversion-exchange relation for the common angle.  Section~\ref{sec:scalar} proves Theorem~\ref{Scalar} by modifying the arguments in \cite{DelPinoWeiYao2015}.
Appendix~\ref{app:TN} contains the detailed spectral and inverse analysis of the exact position matrices and the resonant angle matrices. We write $A_K\asymp B_K$ for two-sided positive
bounds independent of $K$ and $A_K\sim B_K$ for $A_K/B_K\to1$.
All coefficients and potentials are fixed unless specified otherwise.
The self-interaction coefficients are $\mu,\nu$ and the force logarithmic
derivative is always $\nu_*$. The manuscript convention is $A_1$ for the
first cyclic difference and $A_2$ for the second cyclic difference.

\subsection{Use of artificial intelligence}  The idea of adding
$\sum_{l=3}^{N}z_{j,l}\pmb e_l$ at every center $Q_j$ and replacing the action of $SO(2)$ by $SO(N)/SO(N-2)$ is contributed by the human authors.  We then used AI models as research assistants to collect all necessary knowledge for this idea (including the special orthogonal group in $\mathbb{R}^N$ and discrete Green functions for circulant matrices), do calculations and computations, and polish the English.  The final checks remain the sole responsibility of the human authors.

\section{Preliminaries}\label{sec:pre-munu}
In this section, we first recall some well-known facts for the positive solution $w$ of~\eqref{eq:w} and its two relevant weighted eigenvalues.
\begin{lemma}\label{mn:scalarfacts}
(\cite[Lemma~3.1]{DelPinoWeiYao2015})\quad Let $N=2,3$, then the positive solution $w$ of~\eqref{eq:w} is unique up to translations.
It is strictly decreasing in $|x|$ and satisfies the following decay
\begin{eqnarray*}
w(r)=c_N r^{-\frac{N-1}{2}}e^{-r}(1+O(r^{-1}))\quad\text{and}\quad
|\nabla^j w(x)|\lesssim(1+|x|)^{-\frac{N-1}{2}}e^{-|x|},
\end{eqnarray*}
for $1\leq j\leq N$, where $c_N>0$ is a constant depending on $N$.  Moreover, the weighted eigenvalue problem
\begin{equation*}\label{mn:weightedspectrum}
\left\{\aligned
&-\Delta f+f=\lambda w^2f,\\
&f\in H^1(\R^N),
\endaligned\right.
\end{equation*}
has a simple principal eigenvalue $1$, with eigenfunction $w$, no eigenvalue
in $(1,3)$, and eigenspace $\operatorname{span}_{1\leq l\leq N}\{\partial_lw\}$ at $3$, that is, the operator $L_+=-\Delta+1-3w^2$ has Morse index one and
\[
 \ker L_+=\operatorname{span}\{\partial_lw:1\le l\le N\}.
\]
\end{lemma}

The properties of the explicit vector profile and its nonresonant linearization are recalled below.
\begin{proposition}\label{mn:profilekernel}
Let $N=2,3$, $0<\mu\le\nu$ and $\beta\in(0,\mu)\cup(\nu,\infty)$.  Then $(c_{1,\beta}w,c_{2,\beta}w)$, given by \eqref{eqn201}, uniquely solves the limiting system~\eqref{eq:system-lim} and its kernel
consists exactly of the $N$ simultaneous translations, where $c_{1,\beta},c_{2,\beta}$ are defined by~\eqref{eq:amplitudes}.  Moreover, its energy is $(c_{1,\beta}^2+c_{2,\beta}^2)I(w)$.
\end{proposition}
\begin{proof}
The solvability of $(c_{1,\beta}w,c_{2,\beta}w)$ for \eqref{eq:system-lim} follows from \cite[Theorem~1]{Sirakov2007} and \cite[Theorem~1.1]{CLWY26} whereas the kernel assertion follows from Lemma~\ref{mn:scalarfacts} and the direct calculations in \cite[Lemma~3.1]{BartschDancerWang2010}.
\end{proof}

In the resonant case, the situation is slightly different.
\begin{proposition}
Let $N=2,3$.  If $\mu=\nu=\beta>0$, then
the kernel at
$(\mu^{-1/2}w(x)\cos\vartheta, \mu^{-1/2}w(x)\sin\vartheta)$
is
\[
 \operatorname{span}\left\{\mu^{-1/2}e(\vartheta)\partial_lw\text{ and }\mu^{-1/2}Je(\vartheta)w\mid1\le l\le N\right\},
\]
Moreover, for every positive finite-energy solution of the full system, we must have
\begin{equation}\label{mn:compat}
 \int(P(x)-Q(x))uv\, dx=(\mu-\beta)\int u^3v-(\nu-\beta)\int uv^3\, dx.
\end{equation}
\end{proposition}
\begin{proof}
In the resonant case, we observe that the linearized operator $\mathcal L_\beta$, given by \eqref{eq:general-linearization}, is
\[
\mathcal L_\beta=\operatorname{diag}(-\Delta+1-3w^2,-\Delta+1-w^2)
\]
in the orthonormal component basis $(e(\vartheta),Je(\vartheta))$.  Thus, the kernel follows from Lemma~\ref{mn:scalarfacts}. \eqref{mn:compat} follows from testing the first equation by $v$ and the second by $u$ and subtracting.
\end{proof}

For centers $\mathbf Q=(Q_0, Q_1, \cdots, Q_{K-1})$ with
$\rho=\min_{i\ne j}|Q_i-Q_j|\to\infty$, we define
\begin{eqnarray*}\label{norm}
\|(h_1,h_2)\|_{**}=\sup_{x\in\R^N}(\sum_{j=0}^{K-1} e^{-\eta|x-Q_j|})^{-1}(|h_1(x)|+|h_2(x)|)
\end{eqnarray*}
as in \cite{DelPinoWeiYao2015}, where $\eta>0$ is small enough to compare with the positive limiting mass and the
fixed lower bounds for the potentials.  We also periodically extend $\mathbf Q$ to $\mathbb{Z}$ to be a cyclic sequence.
\begin{lemma}\label{mn:interactions}
Let $N=2,3$ and let $\mathbf Q$ be a polygonal configuration given by~\eqref{eq:synccenters}, with perturbations $\overline{\mathbf q}$ bounded in the discrete norm given by \eqref{eq:discretenorm}.  Then
\[\left\{\aligned
&\sup_{x\in\R^N} \sum_{j=0}^{K-1}e^{-\eta |x-Q_j|}\lesssim1,\\
&\int_{\R^N}(\sum_{j=0}^{K-1}e^{-\eta |x-Q_j|})^2\,dx\lesssim K,\\
&\sup_j\sum_{i=1\atop
 i\ne j}^{K}e^{-\eta|Q_i-Q_j|}\lesssim e^{-c\rho}
\endaligned\right.\]
for a universal constant $c>0$.
Moreover, with $H_w$ as in~\eqref{eq:interactions-general} and $\Psi_w=-H_w'$,
\begin{equation}\label{mn:Hder}
 \begin{aligned}
 H_w^{(k)}(d)&=\int_{\R^N} w^3(x)(-\partial_1)^k w(x-d\pmb e_1)\,dx\\
 &=(-1)^k c_Hd^{-(N-1)/2}e^{-d}(1+O(d^{-1}))
 \end{aligned}
\end{equation}
as $d\to\infty$ for $0\leq k\leq 4$ with $c_H=c_N\ds\int_{\R^N} e^{x_1}w^3(x)\, dx>0$ and $c_N$ given by Lemma~\ref{mn:scalarfacts}.  In particular
$$
D=-\frac{d\Psi_w'(d)}{\Psi_w(d)}=d+O(1)
$$
as $d\to\infty$.
\end{lemma}
\begin{proof}
These estimates are given by \cite[Lemma~3.2 and Lemmas~3.4--3.8]{DelPinoWeiYao2015}.
\end{proof}

We next state the projected estimates used in the nonresonant and resonant cases.  We begin with an estimate in the expansion.
First, writing $(u,v)=(U+\phi,V+\psi)$, then solving $S_\beta(u,v)=0$ amounts to solve
\begin{align*}
L(\phi,\psi)+E+N(\phi,\psi)=0,\,\,\,\,\text{where}\,\,\,\,L(\phi,\psi)=S_{\beta}'(U,V)(\phi,\psi)
\end{align*}
and
\begin{align*}
E=S_{\beta}(U,V), N(\phi,\psi)=S_{\beta}(U+\phi,V+\psi)-S_{\beta}(U,V)-S_{\beta}'(U,V)(\phi,\psi)
\end{align*}

\begin{lemma}\label{mn:nonlinear}
Let $N=2,3$ and
\begin{eqnarray*}
S_\beta(u,v)=\left\{\aligned
-\Delta u+P(x)u-\mu u^3-\beta v^2u,\\
-\Delta v+Q(x)v-\nu v^3-\beta u^2v.
\endaligned\right.
\end{eqnarray*}
Then for the pairs $(U,V)$ and $(\phi,\psi)$,
\begin{eqnarray*}
\left\{\aligned
 N_1(\phi,\psi)&=-3\mu U\phi^2-\mu\phi^3-\beta(2V\phi\psi+U\psi^2+\phi\psi^2),\\
 N_2(\phi,\psi)&=-3\nu V\psi^2-\nu\psi^3-\beta(2U\phi\psi+V\phi^2+\psi\phi^2).
\endaligned\right.
\end{eqnarray*}
Moreover, if $\|(U,V)\|_{**}+\|(\phi,\psi)\|_{**}+\|(\widetilde\phi,\widetilde\psi)\|_{**}\lesssim1$, then
\[
\left\{\aligned
&\|N(\phi,\psi)\|_{**}\lesssim\|(\phi,\psi)\|_{**}^2,\\
&\|N(\phi,\psi)-N(\widetilde\phi,\widetilde\psi)\|_{**}
 \lesssim(\|(\phi,\psi)\|_{**}+\|(\widetilde\phi,\widetilde\psi)\|_{**})
          \|(\phi-\widetilde\phi,\psi-\widetilde\psi)\|_{**},
\endaligned\right.
\]
while if $(U,V), (\phi,\psi), (\widetilde\phi,\widetilde\psi)\in\left(H^2(\R^N)\right)^2$ with $\|(U,V)\|_{L^\infty(\R^N)}\lesssim1$, then the analogous estimates hold from $H^2$ to
$L^2$ on a bounded ball.
\end{lemma}
\begin{proof}
The conclusion follows by expanding $(U+\phi)^3$ and $(U+\phi)(V+\psi)^2$ and subtracting their constant and linear terms.
\end{proof}

To state our linear theory in a uniform form, we use $\{X_{j,l}\}$ for the local approximate kernel functions at the ansatz.  More precisely, in the nonresonant case, we denote
\begin{eqnarray*}
X_{j,l}=\binom{\partial_lw_{1,\beta}(x-Q_j)}{\partial_lw_{2,\beta}(x-Q_j)}
\end{eqnarray*}
for $0\leq j\leq K-1$ and $1\leq l\leq N$ while in the resonant case, we denote
\begin{eqnarray}\label{kernel-resonant}
X_{j,l}=\left\{\aligned
&\frac1{\sqrt{\mu}}e(\vartheta_j)\partial_lw_j,\quad 1\le l\le N\text{ and }0\leq j\leq K-1,\\
&\frac1{\sqrt{\mu}}w_jJe(\vartheta_j),\quad l=N+1\text{ and }0\leq j\leq K-1.
\endaligned\right.
\end{eqnarray}
Let
\begin{equation}\label{pf:actualoperator}
S_\beta'(U,V)=
\begin{pmatrix}
-\Delta+P(x)-3\mu U^2-\beta V^2&-2\beta UV\\
-2\beta UV&-\Delta+Q(x)-3\nu V^2-\beta U^2
\end{pmatrix},
\end{equation}
where
\begin{eqnarray}\label{ansatz-nonresonant}
(U,V)=\Bigl(\sum_{j=0}^{K-1}w_{1,\beta}(x-Q_j),\sum_{j=0}^{K-1}w_{2,\beta}(x-Q_j)\Bigr)
\end{eqnarray}
in the nonresonant case and
\begin{eqnarray}\label{ansatz-resonant}
(U,V)=\frac1{\sqrt{\mu}}\Bigl(\sum_{j=0}^{K-1}
 w_j \cos\vartheta_j,\sum_{j=0}^{K-1}
 w_j \sin\vartheta_j\Bigr)
\end{eqnarray}
in the resonant case.  Now, our uniform linear theory can be stated as follows.
\begin{lemma}\label{mn:linearprojection}
Let $N=2,3$.  For the centers $\{Q_j\}$ given by \eqref{eq:synccenters} with $R\asymp K\log K$ for $K\in\mathbb{N}$ sufficiently large, the
projected problem
\begin{equation*}\label{mn:linearprojected}
\left\{\aligned
&S_\beta'(U,V)\binom{\phi}{\psi}=\binom{h}{g}+\sum_{j,l} c_{j,l}X_{j,l},\\
&\int_{\mathbb{R}^N}(\phi,\psi)\cdot X_{j,l}=0,
\endaligned\right.
\end{equation*}
has a unique solution with the estimates
\[
\|(\phi,\psi)\|_{**}+\|{\mathbf c}\|_\infty\lesssim\|(h,g)\|_{**}\quad\text{and}\quad
\|(\phi,\psi)\|_{H^1(\mathbb{R}^N)}\lesssim \|(h,g)\|_{L^2(\mathbb{R}^N)},
\]
where the linearized operator $S_\beta'(U,V)$ is given by \eqref{pf:actualoperator}.
\end{lemma}
\begin{proposition}\label{mn:firstprojection}
Under the assumption of  Lemma \ref{mn:linearprojection}, the nonlinear projected problem
\begin{equation}\label{mn:linearprojected-1}
\left\{\aligned
&S_\beta'(U,V)\binom{\phi}{\psi}=-S_\beta(U,V)-\binom{N_1(\phi,\psi)}{N_2(\phi,\psi)}+\sum_{j,l} c_{j,l}X_{j,l},\\
&\int_{\mathbb{R}^N}(\phi, \psi)\cdot X_{j,l}=0,
\endaligned\right.
\end{equation}
therefore has a
unique small correction $(\phi_K,\psi_K)$  in the nonresonant case
\begin{eqnarray*}
\left\{\aligned
&\|(\phi_K,\psi_K)\|_{**}+\|\mathbf c\|_\infty\lesssim R^{-m_0},\\
&\|(\phi_K,\psi_K)\|_{H^1(\mathbb{R}^N)}\lesssim\sqrt{K}R^{-m_0},\\
&\|(\partial_{\overline{\mathbf q}}\phi_K,\partial_{\overline{\mathbf q}}\psi_K)\|_{**}\lesssim R^{-m_0}.
\endaligned\right.
\end{eqnarray*}
In the resonant case, we have
\begin{eqnarray*}
\left\{\aligned
&\|(\phi_K,\psi_K)\|_{**}+\|\mathbf c\|_\infty\lesssim R^{-m},\\
&\|(\phi_K,\psi_K)\|_{H^1(\mathbb{R}^N)}\lesssim\sqrt{K}R^{-m},\\
&\|(\partial_{\overline{\mathbf q}}\phi_K,\partial_{\overline{\mathbf q}}\psi_K)\|_{**}\lesssim R^{-m},\\
&\|(\partial_{\vartheta}\phi_K,\partial_{\vartheta}\psi_K)\|_{**}+\|(\partial_{\mathbf a}\phi_K,\partial_{\mathbf a}\psi_K)\|_{**}\lesssim R^{1-m}.
\endaligned\right.
\end{eqnarray*}
\end{proposition}
\begin{proof}
The linear estimate follows by a contradiction and blow-up argument (cf. \cite[Lemma~4.3 and Proposition~4.1]{DelPinoWeiYao2015} or \cite[Proposition~3.1]{LiWeiWu2024}).  Thus, we omit it here.  We only point out that in the resonant case, after passing to the limit, the linearized operator $S_\beta'(U,V)$ is $\operatorname{diag}(-\Delta+1-3w^2,-\Delta+1-w^2)$ at every center $Q_j$ in the orthonormal component basis $(e(\vartheta_j),Je(\vartheta_j))$.

The conclusion for the nonlinear projected problem follows from the estimates
$$\|S_\beta(U,V)\|_{**}\lesssim R^{-m_0}$$
 in the nonresonant case and 
$$\|S_\beta(U,V)\|_{**}\lesssim \mu^{-1/2}R^{-m}$$ in the resonant case and Lemma~\ref{mn:nonlinear}.

These estimates can be found in  \cite[Lemma~4.4]{DelPinoWeiYao2015} or \cite[Lemma~2.2]{WangYe2025}, so we also omit them here.  The $C^1$-dependence of the parameters can be proved by adapting the same calculations in the proof of \cite[Proposition~4.1]{DelPinoWeiYao2015}, with the additional local angle function at resonance.
We define
\begin{align*}
B_k=\Bigl\{(\phi,\psi)\in L^\infty(\R^N)\times L^\infty(\R^N)\mid \|(\phi,\psi)\|_{**}\leq C_0 R^{-{m_0}}\Bigr\}.
\end{align*}
Now, we define a map
\begin{align*}
B(\phi,\psi)=-L^{-1}(E+N(\phi,\psi)).
\end{align*}
Using Lemma~\ref{mn:linearprojection}, $B$ is a contraction mapping on $B_k$, Then \eqref{mn:linearprojected-1} holds.
\end{proof}

To use uniform notation, we introduce
\begin{eqnarray*}
\overline{\mathbf p}=\left\{\aligned
&\overline{\mathbf q},\quad\text{in the nonresonant case},\\
&(\overline{\mathbf q},\vartheta,{\mathbf a}),\quad\text{in the resonant case}.
\endaligned\right.
\end{eqnarray*}
\begin{lemma}\label{mn:energyderivative}
Let $(U, V)$ be the ansatz given by \eqref{ansatz-nonresonant} in the nonresonant case and given by  \eqref{ansatz-resonant} in the resonant case.  Then for the corrected projected family $Y=(U+\phi_K, V+\psi_K)$ with
$Y_i=\partial_{p_i}(U, V)$, we have
\begin{equation}\label{mn:gradient}
 \partial_{p_i}I_\beta(Y)
 =\sum_{j,l} c_{j,l}\left(\langle X_{j,l},Y_i\rangle
                         -\langle(\phi, \psi),\partial_{p_i}X_{j,l}\rangle\right).
\end{equation}
In particular,
\begin{equation*}\label{mn:thetaidentity}
 \partial_\vartheta I_\mu(Y)=\int_{\mathbb{R}^N}(Q(x)-P(x))(U+\phi_K)(V+\psi_K) \, dx
\end{equation*}
in the resonant case.
\end{lemma}
\begin{proof}
We denote $\Phi_K=(\phi_K, \psi_K)$ for the sake of simplicity.  Then differentiating $\langle\Phi_K, X_{j,l}\rangle=0$ gives
\[
 \langle\partial_{p_i}\Phi_K,X_{j,l}\rangle=-\langle\Phi_K,\partial_{p_i}X_{j,l}\rangle,
\]
which, together with $S_\beta(Y)=\sum_{j,l}c_{j,l}X_{j,l}$ proved in Proposition~\ref{mn:firstprojection}, implies \eqref{mn:gradient}.

For $\partial_{\vartheta}I_\mu(Y)$, by $\partial_{\vartheta}X_{j,l}=JX_{j,l}$ and $J^T=-J$, we have $\langle\partial_{\vartheta}\Phi_K-J\Phi_K,X_{j,l}\rangle=0$, which, together with $\partial_{\vartheta}(U, V)^T=J(U, V)^T$ and Proposition~\ref{mn:firstprojection}, implies that
\begin{eqnarray*}
\partial_\vartheta I_\mu(Y)&=&\langle S_\mu(Y),JY\rangle\\
&=&\int_{\R^N}((V+\psi)\Delta (U+\phi)-(U+\phi)\Delta (V+\psi))\, dx\\
&&+\int_{\R^N}(Q(x)-P(x))(U+\phi)(V+\psi)\, dx\\
&=&\int_{\R^N}(Q(x)-P(x))(U+\phi)(V+\psi)\, dx.
\end{eqnarray*}
This completes the proof.
\end{proof}

\section{The nonresonant case}\label{sec:nonres-munu}
By Proposition~\ref{mn:firstprojection}, we have constructed a solution of the following projected equation~\eqref{mn:linearprojected-1}.  Thus, the major task in this section is to kill all its multipliers.
\begin{proposition}\label{pf:projection}
Let $N=2,3$. If $\nabla_{\overline{\mathbf q}} I_\beta(Y)=\bf 0$ for some $\|\overline{\mathbf q}\|_*<2$,  then the multipliers in \eqref{mn:linearprojected-1} satisfy $\mathbf c=\mathbf 0$.  In particular, $Y=(U+\phi_K,V+\psi_K)$ is a solution of \eqref{eq:system}.
\end{proposition}
\begin{proof}
By Lemma~\ref{mn:energyderivative},
\begin{equation}\label{pf:gradient-1}
\nabla_{\overline{\mathbf q}}I_\beta(Y)
 =\sum_{j,l}c_{j,l}
 \left(\left\langle X_{j,l},\nabla_{\overline{\mathbf q}}(U,V)\right\rangle
       -\left\langle(\phi_K,\psi_K),\nabla_{\overline{\mathbf q}}X_{j,l}\right\rangle\right).
\end{equation}
Since the first matrix
in~\eqref{pf:gradient-1} is the Gram matrix up to these invertible coordinate
changes and the other matrix has entries $O(\|(\phi_K,\psi_K)\|_{**})$, the conclusion then follows from Proposition~\ref{mn:firstprojection} and \eqref{mn:linearprojected-1}.
\end{proof}

By \eqref{mn:Hder},
\begin{equation*}\label{pf:Hdifferentiate}
\begin{aligned}
H_w^{(k)}(d)&=\int_{\R^N} w^3(x)(-\partial_1)^k w(x-d\pmb e_1)\,dx,\\
&=(-1)^kc_Hd^{-\frac{N-1}{2}}e^{-d}(1+O(d^{-1}))
\end{aligned}
\end{equation*}
for $0\leq k\leq4$ as $d\to+\infty$, which supplies the derivatives needed in addition to the leading
interaction of~\cite[Lemma~3.2]{DelPinoWeiYao2015}.
\begin{lemma}\label{pf:radius}
Let $N=2,3$.  Then in the nonresonant case, the balance
\eqref{eq:balance} has a solution with $d\to+\infty$ under the sign conditions~\eqref{eq:signcondition}, which is unique for $d>0$ sufficiently large with $R$ satisfying~\eqref{eq:R-asymptotic} and
\begin{equation}\label{eq:d-asymptotic}
d_K=m_0\log K+\left(m_0-\frac{N-3}{2}\right)
                   \log(m_0\log K)+O(1).
\end{equation}
Moreover,
\begin{eqnarray*}
e^{-d_K}\asymp R_K^{-m_0}d^{(N-3)/2}, \quad \Psi(d_K)\asymp R_K^{-m_0}/d_K,\quad b_{*,K}\asymp R_K^{-m_0-2},\quad D_K=d_K+O(1)
\end{eqnarray*}
and $\nu_{*,K}\to m_0+1>3$, where $\Psi(d_K)=-(c_{1,\beta}^2+c_{2,\beta}^2)H'_w(d_K)$ and $b_{*,K}, D_K, \nu_{*,K}$ are given by \eqref{eq:D-nueff}.
\end{lemma}
\begin{proof}
The dominant coefficient of $B(r)$ is positive under precisely the
three alternatives in~\eqref{eq:signcondition}, where $B(r)$ is given by \eqref{eq:potential-coefficients}. Thus,
\[
B(R)\asymp R^{-m_0-1}>0\quad\text{and}\quad
\nu_*=-\frac{RB'(R)}{B(R)}\to m_0+1>3.
\]
Since $\Psi(d)\asymp d^{-\frac{N-1}{2}}e^{-d}$ and $2s=d/R$, the balance~\eqref{eq:balance} implies that the solution must satisfy
$$
R_K^{-m_0}\asymp d_K^{-\frac{N-3}{2}}e^{-d_K},
$$
which is equivalent to
$$
e^{-d_K}\asymp R_K^{-m_0}d_K^{\frac{N-3}{2}}.
$$
By taking logarithms on both sides and using $R_K=d_K/(2\sin(\pi/K))$, we have
\begin{align*}
d_K&=m_0\log R_K-\frac{N-3}{2}\log d_K+O(1)\\
 &=m_0\log K+\left(m_0-\frac{N-3}{2}\right)\log d_K+O(1)\\
 &=m_0\log K+\left(m_0-\frac{N-3}{2}\right)\log(m_0\log K)+O(1),
\end{align*}
which, together with dividing by $2\sin(\pi/K)$, implies
\begin{align*}
R_K=\frac{m_0}{2\pi}K\log K
 +\frac{m_0-(N-3)/2}{2\pi}K\log(m_0\log K)+O(K).
\end{align*}

Its logarithmic derivative is
\[
\frac d{dR}\log\frac{2s\Psi(2sR)}{B(R)}
 =\frac{\nu_*-D}{R}<0
\]
for $d>0$ sufficiently large since $D=d+O(1)$, which implies the existence and uniqueness of the solution to the balance~\eqref{eq:balance}.
\end{proof}

Let $\mathcal J({\bf Q})=I_\beta(U+\phi_K,V+\psi_K)$ be the corrected projected energy.
\begin{lemma}\label{mn:force}
Let $N=2,3$.  Then in the nonresonant case,
\begin{equation}\label{pf:fullforce}
-\partial_{Q_j}\mathcal J({\bf Q})
 =B(|Q_j|)\frac{Q_j}{|Q_j|}
  +\sum_{i=j\pm1}\Psi(|Q_i-Q_j|)\frac{Q_i-Q_j}{|Q_i-Q_j|}+r_j
\end{equation}
for all $0\leq j\leq K-1$.  Moreover,
\begin{equation}\label{mn:forceerror}
\|{\bf r}\|_\infty+
 \|{\bf \partial_{\overline{\mathbf q}} r}\|_{\ell^\infty\to\ell^\infty}
 \lesssim R_K^{-2-m_0}\delta_K
\end{equation}
for $\|\overline{\mathbf q}\|_*\le1$.
\end{lemma}
\begin{proof}
Recall that $w_j=w(\cdot-Q_j)$ and $W=\sum_{j=1}^{K}w_j$.  Then
\begin{eqnarray*}\label{mn:scalarenergy}
I_\beta\left(U,V\right)
&=&\frac{c_{1,\beta}^2+c_{2,\beta}^2}{2}
 \int_{\mathbb{R}^N}\left(\left|\nabla W\right|^2+W^2\right)\, dx\notag\\
 &&+\frac12\int_{\mathbb{R}^N}(c_{1,\beta}^2(P(x)-1)+c_{2,\beta}^2(Q(x)-1))W^2\, dx\notag\\
&&-\frac{c_{1,\beta}^2+c_{2,\beta}^2}{4}\int_{\mathbb{R}^N}W^4\, dx.
\end{eqnarray*}
By applying the Taylor expansion and using the tail assumption $(C_2)$, we have
\begin{align*}
\int_{\R^N} (P(x)-1)w_j\partial_{Q_j} w_j\, dx
 &=\frac12\partial_{Q_j}\ds\int_{\R^N}(P(x)-1)w(x-Q_j)^2\,dx\\
 &=-\frac{am\ds\int_{\R^N}w^2}{2}|Q_j|^{-m-2}Q_j
          +O(R_K^{-m-\sigma}+R_K^{-m-3})
\end{align*}
and
\begin{align*}
\int_{\R^N} (Q(x)-1)w_j\partial_{Q_j} w_j\, dx
 &=\frac12\partial_{Q_j}\int_{\R^N}(Q(x)-1)w(x-Q_j)^2\,dx\\
 &=-\frac{bn\ds\int_{\R^N}w^2}{2}|Q_j|^{-n-2}Q_j
          +O(R_K^{-n-\sigma}+R_K^{-n-3}).
\end{align*}
For $i\ne j$,
\[
3\int_{\R^N} w_j^2w_i\partial_{Q_j} w_j\, dx
 =\Psi_w(|Q_i-Q_j|)\frac{Q_i-Q_j}{|Q_i-Q_j|},
\qquad
\]
where $\Psi_w(s)=-H_w'(s)$ and $H_w(s)$ is given by \eqref{mn:Hder}.
Consequently, the common projection is weighted by
$c_{1,\beta}^2+c_{2,\beta}^2$.
The high-order core remainders obey the convolution estimates
of \cite[Lemmas~3.7--3.8 and 5.1]{DelPinoWeiYao2015}.
In particular, for a possibly smaller fixed positive $\eta$,
\begin{eqnarray*}
R_K^{-m_0}\sum_{i\ne j}\int w_iw_j\lesssim R_K^{-m_0}e^{-(1-\eta/4)d}
          \lesssim R_K^{-(2-\eta)m_0},
\end{eqnarray*}
\begin{eqnarray*}
\sum_{i\notin\{j,j-1,j+1\}}
 \bigl(|\Psi(|Q_i-Q_j|)|+|\Psi'(|Q_i-Q_j|)|\bigr)
\lesssim e^{-(2-\eta/4)d_K}
          \lesssim R_K^{-(2-\eta)m_0},
\end{eqnarray*}
since $\sin(l\pi/K)\geq\frac{2l}{\pi}\sin(\pi/K)$
for $1\le l\le K/2$. Fixed powers of $d_K\asymp\log R_K$ are absorbed by the unused fraction
of $\eta$.

It remains to justify replacing the uncorrected projection by the exact
energy derivative. By self-adjointness of $S_\beta'(U,V)$, the single-core kernel identity in \eqref{mn:linearprojected-1}, Lemma~\ref{mn:nonlinear} and Proposition~\ref{mn:firstprojection}, we have
\begin{eqnarray*}
\left\{\aligned
&\left|\left\langle S_\beta'(U,V)\binom{\phi}{\psi},X_{j,l}\right\rangle\right|=|\langle(\phi,\psi),S_\beta'(U,V)X_{j,l}\rangle|
 \lesssim R_K^{-(2-\eta)m_0},\\
&|\langle N((\phi,\psi)),X_{j,l}\rangle|\lesssim R_K^{-2m_0},\\
&\Bigr|\sum_{i,r}c_{i,r}\langle(\phi,\psi),
                  \partial_{Q_{j,l}}X_{i,r}\rangle\Bigl|
 \lesssim R_K^{-2m_0}
\endaligned\right.
 \end{eqnarray*}
for all $0\leq j\leq K-1$ and $1\leq l\leq N$.

Since each differentiated kernel has the same center as the parameter, by Proposition~\ref{mn:firstprojection} and similar to \cite[Lemma~5.1]{DelPinoWeiYao2015}, differentiating the displayed pairings once gives the same orders for the
remaining terms, which proves \eqref{mn:forceerror} and its stated parameter estimates.
\end{proof}

Let $Q_j^0=R\mathcal{A} \pmb n_j$ and ${\bf Q}^0=\left(Q_0^0, Q_1^0, \cdots, Q_{K-1}^0\right)$.  The following lemma provides the exact linearization of the principal force in the total force \eqref{pf:fullforce} under the local frame $(\mathcal A\pmb n_j, \mathcal A\pmb t_j)$ for $N=2$ and $(\mathcal A\pmb n_j, \mathcal A\pmb t_j, \mathcal A\pmb e_3)$ for $N=3$.
\begin{lemma}\label{mn:matrixlinear}
Let $N=2,3$.  Then the principal force of \eqref{pf:fullforce} is given by
\begin{equation*}\label{pf:principalforce}
F_j({\bf Q})=B(|Q_j|)\frac{Q_j}{|Q_j|}
  +\sum_{i=j\pm1}\Psi(|Q_i-Q_j|)\frac{Q_i-Q_j}{|Q_i-Q_j|},
\end{equation*}
which has a zero constant term for $0\leq j\leq K-1$.  Moreover,
\begin{eqnarray*}
\nabla _{\overline{\mathbf q}}F_j\left({\bf Q}^0\right)\overline{\mathbf q}=b_{*,K}(T_N\overline{\mathbf q})_j
\end{eqnarray*}
in the local frames
$(\mathcal A \pmb n_j,\mathcal A \pmb t_j,\mathcal A \pmb e_3)$, with
$T_N$ as in~\eqref{s:eq:T3}--\eqref{eq:matrixcoefficients} and $b_{*,K}$ given by Lemma~\ref{pf:radius}.
\end{lemma}
\begin{proof}
By \eqref{eq:D-nueff} and Lemma~\ref{pf:radius},
\begin{align*}
\nabla_{\overline{\mathbf q}}\left(B(|Q_j|)\frac{Q_j}{|Q_j|}\right)\overline{\mathbf q}=b_{*,K}(-\nu_{*,K}f_j\mathcal A\pmb n_j+g_j\mathcal A\pmb t_j+z_j\mathcal A\pmb e_3)
\end{align*}
and $\frac{\Psi(d_K)}{d_K}=b_{*,K}/(4s^2)$.  On the other hand, the two reference edges at the $j$-th center $Q_j$, in the $j$-th local frame $(\mathcal A\pmb n_j, \mathcal A\pmb t_j)$, are $$\pmb e_+=\left\{\aligned &d_K(-s,c),  \,\, for \,\, N=2\\ & d_K(-s,c,0) \,\, for \,\, N=3\endaligned\right.  \quad\text{and} \quad
\pmb e_-=\left\{\aligned &d_K(-s,-c),  \,\, for \,\, N=2,\\ & d_K(-s,-c,0) \,\, for \,\, N=3,\endaligned\right.$$ where $s=\sin(\pi/K)$, $c=\cos(\pi/K)$ and $d_K$ is given by \eqref{eq:d-asymptotic}.

Moreover, the corresponding changes of the planar frame $(\mathcal A\pmb n_{j+1}, \mathcal A\pmb t_{j+1})$ and $(\mathcal A\pmb n_{j-1}, \mathcal A\pmb t_{j-1})$ into $(\mathcal A\pmb n_j, \mathcal A\pmb t_j)$ are
\begin{eqnarray*}\label{R_+}
R_+=
\begin{pmatrix}c^2-s^2&-2sc\\2sc&c^2-s^2\end{pmatrix}\quad\text{and}\quad R_-=\begin{pmatrix}c^2-s^2&2sc\\-2sc&c^2-s^2\end{pmatrix}\quad\text{for }N=2
\end{eqnarray*}
\begin{eqnarray*}
R_+=\begin{pmatrix}c^2-s^2&-2sc&0\\2sc&c^2-s^2&0\\
0&0&1\end{pmatrix}\quad\text{and}\quad R_-=\begin{pmatrix}c^2-s^2&2sc&0\\-2sc&c^2-s^2&0\\
0&0&1\end{pmatrix}\quad\text{for }N=3.
\end{eqnarray*}
Since
\begin{align*}
\nabla_y\left(\Psi(|y|)\frac {y}{|y|}\right)=\frac{\Psi(|y|)}{|y|}\left(I-(D+1)\frac {y}{|y|}\otimes\left(\frac {y}{|y|}\right)^T\right),
\end{align*}
where $I$ is the identity matrix and $D_K$ is given by \eqref{eq:D-nueff},
by direct calculations and Lemma~\ref{pf:radius},for $N=2$,
\begin{eqnarray*}
\nabla_{\overline{\mathbf q}}F_j({\bf Q})&=&\frac{b_{*,K}}{4s^2}\left(M_+R_+\binom{f_{j+1}}{g_{j+1}}
+
M_-R_-\binom{f_{j-1}}{g_{j-1}}-(M_++M_-)\binom{f_{j}}{g_{j}}\right)\\
&&+b_{*,K}(-\nu_{*,K}f_j\mathcal A\pmb n_j+g_j\mathcal A\pmb t_j)\\
&=&b_{*,K}(F_{j,n}\mathcal A\pmb n_j+F_{j,t}\mathcal A\pmb t_j)
\end{eqnarray*}
 and for $N=3$
\begin{eqnarray*}
\nabla_{\overline{\mathbf q}}F_j({\bf Q})&=&\frac{b_{*,K}}{4s^2}\left(M_+R_+\begin{pmatrix}{f_{j+1}}\\{g_{j+1}}\\{z_{j+1}} \end{pmatrix}
+
M_-R_-\begin{pmatrix}{f_{j-1}}\\{g_{j-1}}\\{z_{j-1}} \end{pmatrix}-(M_++M_-)\begin{pmatrix}{f_{j}}\\{g_{j}}\\{z_{j}} \end{pmatrix}\right)\\
&&+b_{*,K}(-\nu_{*,K}f_j\mathcal A\pmb n_j+g_j\mathcal A\pmb t_j+z_j\mathcal A\pmb e_3)\\
&=&b_{*,K}(F_{j,n}\mathcal A\pmb n_j+F_{j,t}\mathcal A\pmb t_j+F_{j,3}\mathcal A\pmb e_3)
\end{eqnarray*}
in the $j$-th local frame $(\mathcal A\pmb n_j, \mathcal A\pmb t_j, \mathcal A\pmb e_3)$,
where $M_\pm=I-\frac{(D_K+1)}{d_K^2}\pmb e_{\pm}\otimes\pmb e_{\pm}^T $ and
\begin{eqnarray*}
\left\{\aligned
F_{j,n}&=a_{e,K}(A_2\pmb f)_j+(D_K-\nu_{*,K})f_j+b_{e,K}(A_1\pmb g)_j,\\
F_{j,t}&=-b_{e,K}(A_1\pmb f)_j+c_{e,K}(A_2\pmb g)_j,\\
F_{j,3}&=z_j+\frac{(A_2\pmb z)_j}{4s^2},
\endaligned\right.
\end{eqnarray*}
with $a_{e,K}$, $b_{e,K}$, $c_{e,K}$ given by \eqref{eq:matrixcoefficients} with $D_K$ given by Lemma~\ref{pf:radius} and $A_1$, $A_2$ the circulant matrices given by \eqref{eq:A1A2-matrices}.  Finally the constant term in the reference force $F_j({\bf Q})$ is
$(B(R_K)-2s\Psi(d_K))\pmb n_j=0$ by the balance~\eqref{eq:balance}.
\end{proof}

To expand the principal force $\boldsymbol F=(F_1, F_2, \cdots, F_K)$, we also need to compute the remainder
\begin{eqnarray*}
\boldsymbol R(\overline{\mathbf q})=\boldsymbol F({\bf Q})-\boldsymbol F({\bf Q}^0)-\nabla _{\overline{\mathbf q}}F\left({\bf Q}^0\right)\overline{\mathbf q}
\end{eqnarray*}
in its linearization, where we also periodically extend $\boldsymbol F$ to $\mathbb{Z}$ to be a cyclic sequence.
\begin{lemma}\label{mn:geometryerror}
For $\overline{\mathbf q}$ and $\widetilde{\mathbf q}$ satisfying $\|\overline{\mathbf q}\|_*\le1$ and $\|\widetilde{\mathbf q}\|_*\le1$, the remainder $\boldsymbol R(\overline{\mathbf q})$ satisfies
\begin{eqnarray}\label{pf:geobound}
\|\boldsymbol R(\overline{\mathbf q})\|_\infty\lesssim\frac{b_{*,K}d_K^2}{R_K}\|\overline{\mathbf q}\|_*^2\quad\text{and}\quad
\|\boldsymbol R(\overline{\mathbf q})
 -\boldsymbol R(\widetilde{\mathbf q})\|_\infty\lesssim\frac{b_{*,K}d_K^2}{R_K}\|\overline{\mathbf q}-\widetilde{\mathbf q}\|_* .\notag
\end{eqnarray}
\end{lemma}
\begin{proof}
By \eqref{eq:synccenters}, direct calculations show that
\begin{eqnarray*}
Q_{j+1}-Q_j=\left\{\aligned
&(-d_Ks+\omega_{j,n})\mathcal A\pmb n_j+(d_Kc+\omega_{j,t})\mathcal A\pmb t_j,\quad N=2,\\
&(-d_Ks+\omega_{j,n})\mathcal A\pmb n_j+(d_Kc+\omega_{j,t})\mathcal A\pmb t_j+\omega_{j,3}\mathcal A\pmb e_3,\quad N=3,
\endaligned\right.
\end{eqnarray*}
where
\begin{eqnarray*}
\left\{\aligned
\omega_{j,n}&=\cos h\,f_{j+1}-\sin h\,g_{j+1}-f_j,\\
\omega_{j,t}&=\sin h\,f_{j+1}+\cos h\,g_{j+1}-g_j,\\
\omega_{j,3}&=z_{j+1}-z_j
\endaligned\right.
\end{eqnarray*}
with
\begin{equation}\label{eq:wj-estimates}
\max_{0\leq j \leq K-1}|\omega_j|\lesssim h\|\overline{\mathbf q}\|_*\quad\text{and}\quad\max_{0 \leq j \leq K-1}|\omega_{j+1}-\omega_j|\lesssim h^2\|\overline{\mathbf q}\|_*
\end{equation}
and
\begin{eqnarray*}
\omega_j=\left\{\aligned
&\omega_{j,n}\mathcal A\pmb n_j+\omega_{j,t}\mathcal A\pmb t_j,\quad N=2,\\
&\omega_{j,n}\mathcal A\pmb n_j+\omega_{j,t}\mathcal A\pmb t_j+\omega_{j,3}\mathcal A\pmb e_3,\quad N=3.
\endaligned\right.
\end{eqnarray*}
Since
\begin{eqnarray*}
\left|\nabla^j\left(\Psi(|y|)\frac{y}{|y|}\right)\right|\lesssim \Psi(|y|)
\end{eqnarray*}
for any $j\geq1$, by \eqref{mn:Hder}, $\Psi_w(d)=-H_w'(d)$, and $\pmb n_{j+1}=\cos h\pmb n_{j}+\sin h\pmb t_{j}$, we have
\begin{eqnarray*}
&&\sum_{i=j\pm1}\Psi(|Q_i-Q_j|)\frac{Q_i-Q_j}{|Q_i-Q_j|}\\
&=&-2s\Psi(d_K)\mathcal{ A}\pmb n_j+\sum_{i=j\pm1}\nabla_{\overline{\mathbf q}}\left(\Psi(|Q_i^0-Q_j^0|)\frac{Q_i^0-Q_j^0}{|Q_i^0-Q_j^0|}\right)\omega_j+O(\Psi(d_K)h^3\|\overline{\mathbf q}\|_*^2),
\end{eqnarray*}
where $h=\frac{2\pi}{K}$.
Similarly, we also have
\begin{eqnarray*}
B(|Q_j|)\frac{Q_j}{|Q_j|}=B(R_K) \mathcal{A} \pmb n_j+b_{*,K}(-\nu_{*,K} f_j \mathcal{A} \pmb n_j+g_j\mathcal{A} \pmb t_j+z_j\mathcal{A} \pmb e_3)+O\left(b_{*,K}R_K^{-1}\|\overline{\mathbf q}\|_*^2\right).
\end{eqnarray*}
By Lemma~\ref{pf:radius}, summarizing the above estimates gives $\|\boldsymbol R(\overline{\mathbf q})\|_\infty\lesssim\frac{b_{*,K}d_K^2}{R_K}\|\overline{\mathbf q}\|_*^2$.  By Lemma~\ref{mn:interactions}, the Lipschitz estimate can be obtained by a very similar argument, so we omit it here.
\end{proof}

We are now in the position to prove Theorem~\ref{Nonresonant}.

\begin{proof}[Proof of Theorem~\ref{Nonresonant}]
By Lemmas~\ref{pf:projection} and \ref{mn:force}, we only need to construct a critical point of $\mathcal J({\bf Q})$.  We divide this construction into two steps.

\vskip0.06in

{\bf Step.~1}\quad We prove that for each $\mathcal A\in SO(N)$, there is a unique small
$\overline{\mathbf q}(\mathcal A)\in\mathcal E_K$, which is  $C^1$ in $\mathcal A$, satisfying
\begin{equation}\label{pf:nrrangeeq}
\Pi_{\mathcal E_K}\nabla_{{\bf Q}}\mathcal J({\bf Q})=\mathbf 0\quad\text{and}\quad\|\overline{\mathbf q}(\mathcal A)\|_*\lesssim\delta_K.
\end{equation}

Indeed, by Lemmas~\ref{pf:radius}, \ref{mn:force} and~\ref{mn:geometryerror},
\begin{eqnarray}\label{pf:normalizednr}
\left\{\aligned
&\|b_{*,K}^{-1}\nabla_{{\bf Q}}\mathcal J({\bf Q})
             +T_N\overline{\mathbf q}\|_\infty\lesssim\delta_K+\frac{d_K^2}{R_K}\|\overline{\mathbf q}\|_*^2,\\
&\|\partial_{\overline{\mathbf q}}
 (b_{*,K}^{-1}\nabla_{{\bf Q}}\mathcal J({\bf Q})
             +T_N\overline{\mathbf q})\|_{*\to\infty}\lesssim\delta_K+\frac{d_K^2}{R_K}\|\overline{\mathbf q}\|_* ,
\endaligned\right.
\end{eqnarray}
where $\delta_K$ is given by \eqref{eq:deltaK} and $d_K$, $R_K$, $b_{*,K}$ are given by Lemma~\ref{pf:radius}.  It follows from Proposition~\ref{tn:TNestimate} that for each $\mathcal A\in SO(N)$, $\Pi_{\mathcal E_K}\nabla_{{\bf Q}}\mathcal J({\bf Q})=\mathbf 0$ has a solution
\begin{equation*}\label{pf:nrfixedpoint}
\overline{\mathbf q}(\mathcal A)
 =(T_N|_{\mathcal E_K})^{-1}\Pi_{\mathcal E_K}
 \left(b_{*,K}^{-1}\nabla_{\overline{\mathbf q}}\mathcal J({\bf Q})
                         +T_N\overline{\mathbf q}(\mathcal A)\right)\in\mathcal E_K
\end{equation*}
which, together with \eqref{pf:normalizednr}, implies that $\|\overline{\mathbf q}(\mathcal A)\|_*\lesssim\delta_K$.  Since $\overline{\mathbf q}(\mathcal A)\in\mathcal E_K$, by \eqref{pf:normalizednr}, Lemma~\ref{mn:geometryerror}, Proposition~\ref{tn:TNestimate} and the implicit function theorem, $\overline{\mathbf q}(\mathcal A)$ is $C^1$ in $\mathcal A$.

\vskip0.06in

{\bf Step.~2}\quad We prove that there exists $\mathcal A\in SO(N)$ such that $\nabla_{{\bf Q}}\mathcal J({\bf Q})=0$.

Indeed, by Step.~1, we have proved that
\begin{equation}\label{pf:remaininggradient}
\nabla_{{\bf Q}}\mathcal J({\bf Q})
                 =\sum_{a=1}^{\frac{N(N-1)}{2}}\gamma_a\mathbf k_a,
\end{equation}
where $\mathbf k_a$ are the kernel vectors of $T_N$ given by \eqref{ker}.  We denote $F(\mathcal A)=\mathcal J(\mathcal A,\overline{\mathbf q}(\mathcal A))$, then by Step.~1, $F(\mathcal A)$ is $C^1$ in $\mathcal A\in SO(N)$.  Moreover, since $SO(N)$ is compact, $F(\mathcal A)$ has a global minimizer and a global maximizer.  Thus, $F(\mathcal A)$ has at least one critical point $\mathcal A_0$.
By \eqref{eq:synccenters}, $\partial_{\mathbf k_1}{\bf Q}=(-g_j,R+f_j,0)_j$ for $N=2$ and
\begin{eqnarray}\label{pf:rotationvelocities}
\left\{\aligned
&\partial_{\mathbf k_1}{\bf Q}=(-g_j,R+f_j,0)_j,\\
&\partial_{\mathbf k_2}{\bf Q}=(-z_j\sin\theta_j,-z_j\cos\theta_j,
              (R+f_j)\sin\theta_j+g_j\cos\theta_j)_j,\\
&\partial_{\mathbf k_3}{\bf Q}=(z_j\cos\theta_j,-z_j\sin\theta_j,
              -(R+f_j)\cos\theta_j+g_j\sin\theta_j)_j
\endaligned\right.
\end{eqnarray}
for $N=3$.  It follows from $\partial_{\mathbf k_a}\overline{\mathbf q}(\mathcal A)\in\mathcal E_K$ for all $\mathcal A\in SO(N)$,
\begin{eqnarray*}
T_{\mathcal A_0}SO(N)=\text{span}_{1\leq a\leq \frac{N(N-1)}{2}}\{\mathcal A_0\mathbf k_a\}
\end{eqnarray*}
and \eqref{pf:remaininggradient} that
\begin{eqnarray}\label{eqn0003}
\mathbf 0&=&\nabla F(\mathcal A_0)=\partial_{\mathbf k_1}F(\mathcal A_0)
\left(\nabla_{{\bf Q}}\mathcal J({\bf Q})
   \left(\partial_{\mathbf k_a}{\bf Q}
              +\mathcal A_0\partial_{\mathbf k_a}\overline{\mathbf q}(\mathcal A_0)\right)\right)_a\cr
   &=&\gamma_1\,
              \mathbf k_1\cdot\partial_{\mathbf k_1}{\bf Q},
\end{eqnarray}
for $N=2$ and 
\begin{eqnarray*}
\nabla F(\mathcal A_0)&=&(\partial_{\mathbf k_1}F(\mathcal A_0), \partial_{\mathbf k_2}F(\mathcal A_0), \partial_{\mathbf k_3}F(\mathcal A_0))\notag\\
&=&\left(\nabla_{{\bf Q}}\mathcal J({\bf Q})
   \left(\partial_{\mathbf k_a}{\bf Q}
              +\mathcal A_0\partial_{\mathbf k_a}\overline{\mathbf q}(\mathcal A_0)\right)\right)_a\notag\\
&=&
\left(\sum_{b=1}^{3}\gamma_b\,\mathbf k_b\cdot\partial_{\mathbf k_1}{\bf Q}, \sum_{b=1}^{3}\gamma_b\,\mathbf k_b\cdot\partial_{\mathbf k_2}{\bf Q}, \sum_{b=1}^{3}\gamma_b\,\mathbf k_b\cdot\partial_{\mathbf k_3}{\bf Q}\right)
\end{eqnarray*}
for $N=3$.
Since by \eqref{pf:nrrangeeq} and \eqref{orth-ker}, the matrix ${\mathbf M}=\left(\mathbf k_b\cdot\partial_{\mathbf k_a}{\bf Q}\right)$ is diagonally dominant, by \eqref{eqn0003}, we must have
\begin{eqnarray*}
\left\{\aligned
&\gamma_1=0,\quad N=2,\\
&(\gamma_1, \gamma_2, \gamma_3)=\mathbf 0,\quad N=3,
\endaligned\right.
\end{eqnarray*}
which, together with \eqref{pf:remaininggradient}, implies that $\nabla_{{\bf Q}}\mathcal J({\bf Q})=0$ for $\mathcal A_0$.
\end{proof}

\section{The resonance case}\label{sec:res-munu}
In this section, we consider the resonance case $\mu=\nu=\beta>0$, under $(C_2)$ with $m=n$ and $a=b>0$.  The profile is given by
\eqref{ansatz-resonant}.  By Proposition~\ref{mn:firstprojection} we construct
a solution of the projected equation~\eqref{mn:linearprojected-1}
\begin{lemma}
Let $N=2,3$ and $\mu=\nu=\beta$. If $\partial_{\overline{\mathbf q}} I_\beta(Y)=\bf 0$ and $\partial_{\vartheta_j} I_\beta(Y)=0$ for some $\|\overline{\mathbf q}\|_*<2$ and all $0\leq j\leq K-1$, where $\|\cdot\|_*$ is the discrete norm for cyclic sequence given by \eqref{eq:discretenorm}, then the multipliers in \eqref{mn:linearprojected-1} satisfy $\mathbf c=\mathbf 0$.  In particular, $Y(\overline{\mathbf q})=(U+\phi_K,V+\psi_K)$ is a solution of \eqref{eq:system}.
\end{lemma}
\begin{proof}
By Lemma~\ref{mn:energyderivative} and $\partial_{\vartheta}(U, V)=J(U, V)$,
\begin{equation}\label{pf:gradient-2}
\nabla_{\overline{\mathbf q}}I_\beta(Y)
 =\sum_{j,l}c_{j,l}
 \left(\left\langle X_{j,l},\nabla_{\overline{\mathbf q}}(U,V)\right\rangle
       -\left\langle(\phi_K,\psi_K),\nabla_{\overline{\mathbf q}}X_{j,l}\right\rangle\right).
\end{equation}
and
\begin{equation}\label{pf:gradient-3}
\partial_{\vartheta_j}I_\beta(Y)
 =\sum_{j,l}c_{j,l}
 \left(\left\langle X_{j,l},J(U, V)\right\rangle
       -\left\langle(\phi_K,\psi_K),\partial_{\vartheta_j}X_{j,l}\right\rangle\right).
\end{equation}
By \eqref{kernel-resonant} and \eqref{ansatz-resonant},
the combining matrix of first matrices
in~\eqref{pf:gradient-2} and \eqref{pf:gradient-3} is the Gram matrix up to these invertible coordinate
changes.  Since the other matrix has entries $O(\|(\phi_K,\psi_K)\|_{**})$, the conclusion then follows from Proposition~\ref{mn:firstprojection} and \eqref{mn:linearprojected-1}.
\end{proof}

The centers ${\bf Q}$ and the
configuration $\overline{\mathbf q}$ in the resonant case are as in
\eqref{eq:synccenters}--\eqref{eq:configuration}.  For the additional rotations on $(U, V)$, we choose a fixed $\vartheta_*\in(0,\pi/4)$ and keep
$\vartheta\in[\vartheta_*,\pi/2-\vartheta_*]$ for all $0\leq j\leq K-1$.
According to  $(C_2)$ with $m=n$ and $a=b>0$ and the new profile in \eqref{ansatz-resonant}, we define
\[
\mathcal J\left({\mathbf Q},\vartheta_1,\cdots,\vartheta_K\right):=I_\mu\left(Y({\mathbf Q},\vartheta_1,\cdots,\vartheta_K)\right)\quad\text{with}\quad Y=(U+\phi_K,V+\psi_K).
\]
\begin{lemma}\label{mn:resforces}
Let $N=2,3$.  Then in the resonant case,
\begin{eqnarray*}\label{pf:resforcefull}
-\partial_{Q_j}\mathcal J=A_0m\frac{Q_j}{|Q_j|^{m+2}}
 +\frac1\mu\sum_{i=j\pm1}
  \Psi_w(|Q_i-Q_j|)\frac{Q_i-Q_j}{|Q_i-Q_j|}
                       \cos(a_i-a_j)+r_{p,j},
\end{eqnarray*}
\begin{eqnarray*}
\partial_{\vartheta_j}\mathcal J=\frac{q_j^*-p_j^*}{2\mu}\sin(2\vartheta_j)
 +\frac1\mu\sum_{i=j\pm1}H_w(|Q_i-Q_j|)
                       \sin(a_j-a_i)+r_{{\bf a},j},
\end{eqnarray*}
where $A_0=\frac{a\int_{\R^N}w^2}{2\mu}$, $\Psi_w(s)=-H_w'(s)$ is given in Lemma~\ref{mn:interactions}, \begin{eqnarray*}
p_j^*=\int_{\R^N}(P(x)-1)w_j^2\, dx\quad\text{and}\quad
q_j^*=\int_{\R^N}(Q(x)-1)w_j^2\, dx.
\end{eqnarray*}
Moreover,
\begin{equation}\label{pf:resanalyticerror}
\left\{\aligned
&\|{ {\bf r}_{ p}}\|_\infty+
 \| {\partial_{\overline{\mathbf q}} {\bf r}_{p}}\|_\infty+\max_{0 \leq j \leq K-1}\| {\partial_{\vartheta_j} {\bf r}_{ p}}\|_\infty\lesssim R_K^{-2-m_0}\delta_K,\\
&\|{{\bf r_{ a}}}\|_\infty+
\|{\partial_{\overline{\mathbf q}} {\bf r_{ a}}}\|_\infty+\max_{0 \leq j \leq K-1}\| \partial_{\vartheta_j} {\bf {r}_{a}}\|_\infty \lesssim R_K^{-2-m_0}\delta_K.
 \endaligned\right.
\end{equation}
\end{lemma}
\begin{proof}
The proof is similar to that of Lemma~\ref{mn:force}, so we only sketch it.  The energy~\eqref{eq:energy} in the resonant case is given by
\begin{equation*}
 I_\mu(u,v)=\frac12\int_{\R^N}\bigl(|\nabla u|^2+|\nabla v|^2+P(x)u^2+Q(x)v^2\bigr)\, dx-\frac\mu4\int_{\R^N}(u^2+v^2)^2\, dx.
\end{equation*}
For the profile $(U,V)$ in the resonant case given by \eqref{ansatz-resonant}, we write
\begin{equation*}\label{eq:energy-Z}
 \begin{aligned}
I_\mu(U,V)={}&\frac1{2}\int_{\R^N}(|\nabla U|^2+|U|^2+|\nabla V|^2+|V|^2)\, dx-\frac{\mu}{4}\int_{\R^N}(U^2+V^2)^2\, dx\\
&+\frac1{2}\int_{\R^N}\bigl((P(x)-1)U^2+(Q(x)-1)V^2\bigr)\, dx.
 \end{aligned}
\end{equation*}
Since for a single profile,
\[
\left|\frac1{\sqrt\mu}w_je(\vartheta_j)\right|^2=\frac1\mu w_j^2
\quad\text{and}\quad e(\vartheta_i)\cdot e(\vartheta_j)
       =\cos(\vartheta_i-\vartheta_j),
\]
by the scalar equation $-\Delta w_i+w_i=w_i^3$ in $\R^N$ and \eqref{mn:Hder}, we can estimate $I_\mu(U,V)$ similarly to the proof of Lemma~\ref{mn:force} to obtain that
\begin{eqnarray*}
I_\mu(U,V)&=&\frac K\mu I(w)+
\frac1{2\mu}\sum_{j=1}^{K}(p_j^*\cos^2\vartheta_j+q_j^*\sin^2\vartheta_j)\\
&&-\frac1\mu\sum_{j=1}^{K}H_w(|Q_{j+1}-Q_j|)
                  \cos(\vartheta_{j+1}-\vartheta_j)+\mathcal R_0,
\end{eqnarray*}
where 
\begin{eqnarray*}
I(w)=\frac{1}{2}\int_{\R^N}|\nabla w|^2+|w|^2\,dx-\frac{1}{4}\int_{\R^N}|w|^4\,dx
\end{eqnarray*}
is the energy of the profile $w$ and $\mathcal R_0$ is the total remainder in the above expansion given by
\begin{align*}\label{eq:R0}
 \mathcal R_0={}&\frac1\mu\sum_{i<j}\int
 \bigl((P(x)-1)\cos\vartheta_i\cos\vartheta_j
 +(Q(x)-1)\sin\vartheta_i\sin\vartheta_j\bigr)w_iw_j\,dx\notag\\
 &-\frac1{4\mu}\int\mathcal P_4(x)\,dx
 -\frac1\mu\sum_{\substack{i<j\\i\not\equiv j\pm1\ ({\rm mod}\ K)}}
 H_w(|Q_i-Q_j|)\cos(\vartheta_i-\vartheta_j)
\end{align*}
with $\mathcal P_4(x)$ the remainder in expanding $\mu^2(U^2+V^2)^2$ given by
\begin{eqnarray*}\label{eq:P4}
 \mathcal P_4&=&\sum_{i<j}\bigl(2+4(e_i\cdot e_j)^2\bigr)w_i^2w_j^2+\sum_i\sum_{\substack{j<k\\j,k\ne i}}
 \bigl(4(e_j\cdot e_k)+8(e_i\cdot e_j)(e_i\cdot e_k)\bigr)
 w_i^2w_jw_k\notag\\
 &&+8\sum_{i<j<k<\ell}
 \Bigl((e_i\cdot e_j)(e_k\cdot e_\ell)
       +(e_i\cdot e_k)(e_j\cdot e_\ell)+(e_i\cdot e_\ell)(e_j\cdot e_k)\Bigr)
 w_iw_jw_kw_\ell.\notag
\end{eqnarray*}
Here, for the sake of simplicity, we denote $e(\vartheta_i)$ by $e_i$.  Since
\begin{eqnarray*}
|e_i|+|\partial_{\vartheta_i} e_i|+|\partial_{\vartheta_i}^2 e_i|\lesssim1,
\end{eqnarray*}
and similar to \eqref{pf:gradient-2} and \eqref{pf:gradient-3},
$\partial_{Q_j}\mathcal R_0$ and $\partial_{\vartheta_j}\mathcal R_0$ share the same decay order of $\mathcal R_0$, the rest of the estimates for $r_{p,j}$ and $r_{{\bf a},j}$ are the same as those for \eqref{mn:forceerror}, the only difference is that we have $m_0=m$ in the resonant case.  Thus, we omit it here.
\end{proof}

We still use the balance
\begin{equation}\label{eq:balance-1}
 \frac{A_0m}{R_K^{m+1}}=\frac{2s}{\mu}\Psi_w(d_K)\quad\text{for}\quad d_K=2R_Ks
\end{equation}
in the resonant case to kill the constant term in $\partial_{\bf Q}\mathcal{J}$, where we recall that $s=\sin\frac{h}{2}$ and $c=\cos \frac{h}{2}$ with $h=\frac{2\pi}{K}$ and $\Psi_w(d_K)$ is given by Lemma~\ref{mn:interactions} with $d_K$ given by Lemma~\ref{pf:radius}.  Thus, the principal force in
$
(\partial_{\bf Q}\mathcal{J}, \partial_{\vartheta_1}\mathcal{J},\cdots, \partial_{\vartheta_K}\mathcal{J})
$
is
\begin{eqnarray*}\label{pf:resforceprinciple}
\left\{\aligned
&F_j({\bf Q},\mathbf a)=A_0m|Q_j|^{-m-2}Q_j
 +\frac1\mu\sum_{i=j\pm1}
  \Psi_w(|Q_i-Q_j|)\frac{Q_i-Q_j}{|Q_i-Q_j|}
                       \cos(a_i-a_j),\\
&G_j({\bf Q},\mathbf a)=\frac{q_j^*-p_j^*}{2\mu}\sin(2(\vartheta+a_j))
 +\frac1\mu\sum_{i=j\pm1}H_w(|Q_i-Q_j|)
                       \sin(a_j-a_i),
                       \label{pf:resangleprinciple}
\endaligned\right.
\end{eqnarray*}
The linearization of $(\boldsymbol F({\bf Q},\mathbf a), \boldsymbol G({\bf Q},\mathbf a))$ at $\overline{\mathbf q}=0$ and $\mathbf a=\mathbf 0$ is given by
\begin{eqnarray}\label{eqn0005}
(\nabla\boldsymbol F({\bf Q},\mathbf a), \nabla\boldsymbol G({\bf Q},\mathbf a))=\text{diag}\left(T_N, -\frac{H_w(d_K)}\mu A_2+\frac{q_j^*-p_j^*}{\mu}\cos2\vartheta\right),
\end{eqnarray}
where $H_w(d_K)$ is given by Lemma~\ref{mn:interactions} and $A_2$ is the circulant matrix given by \eqref{eq:A1A2-matrices}.
\begin{lemma}
For $\|\overline{\mathbf q}\|_*+\|\mathbf a\|_*\le1$, we have
\begin{eqnarray}
\left\{\aligned
&\left\|-b_{*,K}^{-1}\partial_{{\bf Q}}\mathcal J({\bf Q},\mathbf a)
                        -T_N\overline{\mathbf q}\right\|_\infty\lesssim\delta_K+\frac{d_K^2}{R_K}
           (\|\overline{\mathbf q}\|_*^2+\|\mathbf a\|_*^2),\label{pf:positionremainder}\\
&\left\|\frac\mu{4s^2H_w(d_K)}\partial_{\mathbf a}\mathcal J({\bf Q},\mathbf a)
                       +\frac{A_2}{4s^2}\mathbf a\right\|_\infty\lesssim\frac{\delta_K}{d_K}
       +h\|\overline{\mathbf q}\|_*\|\mathbf a\|_*
       +h^2\|\mathbf a\|_*^3,\label{pf:angleremainder}
\endaligned\right.
\end{eqnarray}
where $b_{*,K}$ is given by Lemma~\ref{pf:radius}.
\end{lemma}
\begin{proof}
Since
\begin{eqnarray*}
|1-\cos(a_{j+1}-a_j)|\lesssim h^2\|\mathbf a\|_*^2,\,\,
|\cos(a_{j+1}-a_j)-\cos(a_j-a_{j-1})|\lesssim h^3\|\mathbf a\|_*^2,
\end{eqnarray*}
by \eqref{pf:resanalyticerror} and Lemma~\ref{pf:radius}, the first estimate in \eqref{pf:angleremainder} for
$$
\left\|-b_{*,K}^{-1}\partial_{{\bf Q}}\mathcal J({\bf Q},\mathbf a)-T_N\overline{\mathbf q}\right\|_\infty
$$
is the same as that in Lemma~\ref{mn:geometryerror}, so we omit it here.

For
$$
\left\|\frac\mu{4s^2H_w(d_K)}\partial_{\mathbf a}\mathcal J({\bf Q},\mathbf a)+\frac{A_2}{4s^2}\mathbf a\right\|_\infty,
$$
by \eqref{eq:wj-estimates} and Lemma~\ref{mn:interactions},
\begin{eqnarray*}
&&|H_w(|Q_{j+1}-Q_j|)-H_w(d_K)|
                   |a_{j+1}-2a_j+a_{j-1}|\\
&&+|H_w(|Q_{j+1}-Q_j|)-H_w(|Q_j-Q_{j-1}|)|
                   |a_j-a_{j-1}|\\
&\lesssim&H_w(d_K)h^3\|\overline{\mathbf q}\|_*\|\mathbf a\|_*.
\end{eqnarray*}
Besides, since
\[
|\sin t-t|\lesssim|t|^3\quad\text{and}\quad
|\sin t-t-\sin s+s|\lesssim(t^2+s^2)|t-s|,
\]
by \eqref{eq:wj-estimates} and the balance~\eqref{eq:balance-1}, we can get the second estimate in \eqref{pf:angleremainder}.
\end{proof}

For a cyclic sequence, we define
\[
\Pi_0\mathbf a=\mathbf a-\frac1K\left(\sum_{j=0}^{K-1}a_j\right)\mathbf1.
\]
\begin{proposition}\label{mn:resrange}
For every $\mathcal A\in SO(N)$ and
$\vartheta\in[\vartheta_*,\pi/2-\vartheta_*]$, there is a unique
pair
$\overline{\mathbf q}(\mathcal A,\vartheta)\in\mathcal E_K,
\quad\text{and}\quad \mathbf a(\mathcal A,\vartheta)\in\mathbf1^\perp,
$
satisfying
\begin{equation*}\label{pf:jointbounds}
\|\overline{\mathbf q}\|_*\lesssim\delta_K
\quad\text{and}\quad \|\mathbf a\|_*\lesssim\delta_K/d_K,
\end{equation*}
such that
\begin{equation*}\label{pf:jointrange}
\Pi_{\mathcal E_K}\partial_{{\mathbf Q}}\mathcal J=\mathbf 0,\,\,\,
\Pi_0(\partial_{\vartheta_1}\mathcal J, \partial_{\vartheta_2}\mathcal J, \cdots, \partial_{\vartheta_K}\mathcal J)=\mathbf 0.
\end{equation*}
Moreover, $\overline{\mathbf q}(\mathcal A,\vartheta)$ and $\mathbf a(\mathcal A,\vartheta)$ are $C^1$ in $(\mathcal A,\vartheta)$.
\end{proposition}
\begin{proof}
Since $\frac{q_j^*-p_j^*}{\mu H_w(d_K)}\cos2\vartheta\to0$ as $K\to\infty$ by the balance~\eqref{eq:balance-1}, where $p_j^*$ and $q_j^*$ are given by Lemma~\ref{mn:resforces}, by \eqref{eqn0005} and Lemma~\ref{tn:algebra}, the proof is the same as that for Step.~1 in the proof of Theorem~\ref{Nonresonant}, so we omit it here.
\end{proof}

We are now in the position to prove Theorem~\ref{Resonant}.

\begin{proof}[Proof of Theorem~\ref{Resonant}]
We divide the proof into two steps.

\vskip0.06in

{\bf Step.~1}\quad Let
\begin{equation*}\label{pf:Fres}
F(\mathcal A,\vartheta)
 =\mathcal J(\mathcal A,
       \overline{\mathbf q}(\mathcal A,\vartheta),
       \vartheta,\mathbf a(\mathcal A,\vartheta)),
\end{equation*}
then there exists $\mathcal A_0\in SO(N)$ such that $(\mathcal A_0, \pi/4)$ is a critical point of
$F(\mathcal A, \vartheta)$.

Fix a sufficiently large even $K$, by $P(x)=Q(-x)$, the transformation $\mathscr S(u,v)(x)=(v(-x),u(-x))$ obeys
\[
S_\mu(\mathscr S(u,v))=\mathscr SS_\mu(u,v)\quad\text{and}\quad
I_\mu(\mathscr S(u,v))=I_\mu(u,v).
\]
Moreover, let $\bar j=j+K/2$ modulo $K$. Then the induced parameter transformation $\mathscr S(\overline{\mathbf q}, \vartheta, \mathbf a)=(\overline{\mathbf q}', \vartheta', \mathbf a')$ is given by
\begin{equation}\label{pf:relabel}
\left\{\aligned
&f'_j=f_{\bar j},\quad g'_j=g_{\bar j},\quad
z'_j=-z_{\bar j},\\
&\vartheta'=\pi/2-\vartheta,\quad a'_j=-a_{\bar j},
\endaligned\right.
\end{equation}
with $z$ omitted for $N=2$ and $\mathcal A$ unchanged. Since
$\pmb n_{\bar j}=-\pmb n_j$ and $\pmb t_{\bar j}=-\pmb t_j$,
\[
Q'_j=-Q_{\bar j}\quad\text{and}\quad
\vartheta'+a'_j=\pi/2-(\vartheta+a_{\bar j}).
\]
The uniqueness in Proposition~\ref{mn:firstprojection} gives
$(\phi', \psi')=\mathscr S(\phi, \psi)$ and
\[
\mathcal J(\mathcal A,\overline{\mathbf q}',\pi/2-\vartheta,\mathbf a')
 =\mathcal J(\mathcal A,\overline{\mathbf q},\vartheta,\mathbf a).
\]
By \eqref{pf:relabel}, the transformation $\mathscr S$ preserves all the range constraints
\begin{align*}
&\sum_{j=0}^{K-1}g'_j=\sum_{j=0}^{K-1}g_j=0,\,\,\,\,\sum_{j=0}^{K-1}a'_j=-\sum_{j=0}^{K-1}a_j=0,\\
&\sum_{j=0}^{K-1}z'_j\sin\theta_j=\sum_{j=0}^{K-1}z_j\sin\theta_j=0,\,\,\,\,
\sum_{j=0}^{K-1}z'_j\cos\theta_j=\sum_{j=0}^{K-1}z_j\cos\theta_j=0,
\end{align*}
and the used discrete norms.
The uniqueness in Proposition~\ref{mn:resrange} therefore implies that
\begin{equation*}\label{pf:Fsymmetry}
F(\mathcal A,\vartheta)=F(\mathcal A,\pi/2-\vartheta).
\end{equation*}
Differentiation at
$\vartheta=\pi/4$ gives
$F_\vartheta(\mathcal A,\pi/4)=-F_\vartheta(\mathcal A,\pi/4)=0$.

As in the nonresonant case, since $SO(N)$ is compact, $F(\mathcal A, \pi/4)$ has a global minimizer and a global maximizer.  Thus, by Proposition~\ref{mn:resrange}, $F(\mathcal A, \pi/4)$ has at least one critical point $\mathcal A_0$.  It follows that $(\mathcal A_0, \pi/4)$ is a critical point of
$F(\mathcal A, \vartheta)$.

\vskip0.06in

{\bf Step.~2}\quad We prove that $\partial_{\mathbf Q}\mathcal{J}=\mathbf 0$ and
\begin{eqnarray*}
(\partial_{\vartheta_1}\mathcal J, \partial_{\vartheta_2}\mathcal J, \cdots, \partial_{\vartheta_K}\mathcal J)=\mathbf 0
\end{eqnarray*}
for $\overline{\mathbf q}(\mathcal A_0,\pi/4)$ and $\mathbf a(\mathcal A_0,\pi/4)$.
By Proposition~\ref{mn:resrange},
\begin{eqnarray*}
\partial_{\bf Q}\mathcal{J}=\sum_{b=1}^{\frac{N(N-1)}{2}}\gamma_b\mathbf k_b\quad\text{and}\quad(\partial_{\vartheta_1}\mathcal J, \partial_{\vartheta_2}\mathcal J, \cdots, \partial_{\vartheta_K}\mathcal J)=\varrho \mathbf 1.
\end{eqnarray*}
Since $\sum_ja_j=0$ and $\overline{\mathbf q}(\mathcal A,\vartheta)\in\mathcal{E}_K$, by the computations in \eqref{pf:rotationvelocities},
\begin{eqnarray*}
\mathbf 0&=&\partial_{\mathcal A} F(\mathcal A_0, \pi/4)
=\partial_{\bf Q}\mathcal{J}(\partial_{\mathcal A}{\bf Q}+\mathcal A_0\partial_{\mathcal A}\overline{\mathbf q})+\partial_{\mathbf a}\mathcal{J}\partial_{\mathcal A}\mathbf a
=\partial_{\bf Q}\mathcal{J}\partial_{\mathcal A}{\bf Q}\\
&=&
\left\{\aligned
&\gamma_1\mathbf k_1\cdot\partial_{\mathbf k_1}{\bf Q},\quad N=2,\\
&\left(\sum_{b=1}^{3}\gamma_b\,
              \mathbf k_b\cdot\partial_{\mathbf k_1}{\bf Q}, \sum_{b=1}^{3}\gamma_b\,
              \mathbf k_b\cdot\partial_{\mathbf k_2}{\bf Q}, \sum_{b=1}^{3}\gamma_b\,
              \mathbf k_b\cdot\partial_{\mathbf k_3}{\bf Q}\right),\quad N=3,
\endaligned\right.
\end{eqnarray*}
which implies that $(\gamma_1, \gamma_2, \gamma_3)=\mathbf 0$.  Besides,
\begin{eqnarray*}
0=\partial_{\vartheta} F(\mathcal A_0, \pi/4)
=\partial_{\bf Q}\mathcal{J}\mathcal A_0\partial_{\vartheta}\overline{\mathbf q}+\sum_{j=1}^K\partial_{\vartheta_j}\mathcal{J}+\partial_{\mathbf a}\mathcal{J}\partial_{\vartheta}\mathbf a
=\sum_{j=1}^K\partial_{\vartheta_j}\mathcal{J}
=K\varrho,
\end{eqnarray*}
which implies that $\varrho=0$.  Finally, uniqueness
in the construction also yields $v_K(x)=u_K(-x)$.
\end{proof}

\section{The Scalar version and prove Theorem~\ref{Scalar}}\label{sec:scalar}
In this section, we study \eqref{eq:scalar} and prove Theorem~\ref{Scalar}.

\begin{proof}[Proof of Theorem~\ref{Scalar}]
We begin with the study of the first projected problem
\begin{eqnarray*}
\left\{\aligned
&-\Delta\phi+V(x)\phi-\mu pW^{p-1}\phi=-E-N(\phi)+\sum_{j,l}c_{j,l}\partial_{l}w_j,\\
&\int_{\R^N}\phi\partial_lw_j\,dx=0,\quad 1\leq l\leq N\text{ and }1\leq j\leq K,
\endaligned\right.
\end{eqnarray*}
where $E=(V(x)-1)W -\mu(W^p-\sum_{j=0}^{K-1}w_j^p)$ and $N(\phi)=-\mu((W+\phi)^{p}-W^p-pW^{p-1}\phi)$.  This can be completed by following the arguments and calculations in \cite[Section~4]{DelPinoWeiYao2015} to get the same results, including the existence, uniqueness and decay estimates.  Let
\begin{eqnarray*}
\mathcal H({\bf Q}_{s})=H(W+\phi).
\end{eqnarray*}
We next expand the energy as in Lemma~\ref{mn:force} to get
\begin{eqnarray*}\label{pf:fullforce-scalar}
-\nabla_{Q_{j,s}}\mathcal H({\bf Q}_{s})
=\frac{A_0m_0Q_{j,s}}{|Q_{j,s}|^{m_0+2}}
 +\frac1\mu\sum_{i=j\pm1}
  \Psi_{w,s}(|Q_{i,s}-Q_{j,s}|)\frac{Q_{i,s}-Q_{j,s}}{|Q_{i,s}-Q_{j,s}|}
+r_{j,s}
\end{eqnarray*}
for all $0\leq j\leq K-1$, where $A_0=\frac{a}{2\mu} \int_{\mathbb{R}^N} w^2$,
\begin{eqnarray*}
\Psi_{w,s}(d)=-H_{w,s}'(d)=\int_{\R^N} w^p(x)\partial_1 w(x-d\pmb e_1)\,dx
\end{eqnarray*}
and
\begin{equation*}
\|{\bf r}_{s}\|_\infty+\|{\bf \partial_{\overline{\mathbf q}} r}_{s}\|_\infty\lesssim R_K^{-2-m_0}\delta_K
\end{equation*}
with the balance $ \frac{A_0m_0}{R_K^{m_0+1}}=\frac{2s}{\mu}\Psi_{w,s}(d_K)$ for $d_K=2R_Ks$ and $s=\sin\frac{\pi}{K}$.
By similar calculations for Lemma~\ref{mn:matrixlinear}, we have
\begin{eqnarray*}
\nabla_{\overline{\mathbf q}_{s}} \boldsymbol F=\text{diag}(T_\parallel, \underbrace{T_{\perp}, \cdots, T_{\perp}}_{N-2}),
\end{eqnarray*}
where $\boldsymbol F=(F_0, F_1, \cdots, F_{K-1})$ with
\begin{eqnarray*}
F_j=\frac{A_0m_0Q_{j,s}}{|Q_{j,s}|^{m_0+2}}
 +\frac1\mu\sum_{i=j\pm1}\Psi_{w,s}(|Q_{i,s}-Q_{j,s}|)\frac{Q_{i,s}-Q_{j,s}}{|Q_{i,s}-Q_{j,s}|},
\end{eqnarray*}
Thus, we can solve the second projected problem $\Pi_{\mathcal E_{K,s}}\nabla_{\mathbf Q_{s}}\mathcal H({\bf Q}_{s})=\mathbf 0$ by similar arguments for Step.~1 in the proof of Theorem~\ref{Nonresonant} under the orthogonal condition~\eqref{eq:spatialconditions-scalar} for the configuration $\overline{\mathbf q}_{s}$, after the Lipschitz estimate as that in Lemma~\ref{mn:geometryerror}, where $\mathcal E_{K,s}=(\ker \nabla_{\overline{\mathbf q}_{s}} \boldsymbol F)^\perp$.  Finally, by the same calculations for \eqref{pf:rotationvelocities} and \eqref{eqn0003}, we can kill all multipliers in the first projected problem by choosing a critical point of $F(\mathcal{A})=\mathcal H({\bf Q}_{s}(\mathcal{A}))$.
\end{proof}

\section{Acknowledgements}
The research is supported by NSFC (No. 12161093, No. 12126356, No. 12631008, No. 12671144), Modern Applied Mathematics and Life Sciences Interdisciplinary Research Team in Yunnan (No. 202405AS350003), Yunnan Fundamental Research Projects (No. 202601CJ070001) and Yunnan Revitalization Talent Support Program.

\appendix
\section{General theory of circulant matrices}
\label{app:TN}
In this section, we provide a detailed analysis on the matrix~$T_N$ given by \eqref{s:eq:T3}.

\vskip0.12in

\subsection{Notations}\quad Let $K\ge5$.  We recall that $h=2\pi/K$, $s=\sin(\pi/K)$ and $c=\cos(\pi/K)$.  For a cyclic sequence $\mathbf a=(a_0, a_1, \cdots, a_{K-1})$, its centered first and second differences are
\begin{equation}\label{tn:differences}
 (A_1\mathbf a)_j=a_{j+1}-a_{j-1}\quad\text{and}\quad
 (A_2\mathbf a)_j=a_{j+1}-2a_j+a_{j-1}.
\end{equation}
The discrete norm of $\mathbf a$ remains
\[
 \|\mathbf a\|_* =\|\mathbf a\|_\infty+
 \left\|\left(\frac{a_{j+1}-a_j}{h}\right)_j\right\|_\infty+
 \left\|\frac{A_2\mathbf a}{h^2}\right\|_\infty,
\]
while the Euclidean norm of $\mathbf a$ is denoted by $\|\mathbf a\|_2$.
Block maximum norms of matrices take the maximum over the component blocks.  Clearly, the usage of the
sum of those norms changes constants by at most the fixed number of blocks.  In particular, no factor depending on the parameter $K$ is introduced by this convention.  We also introduce
\[
\bar a=\frac1K\sum_{j=0}^{K-1}a_j
\]
as the average of $\mathbf a$.

\vskip0.12in

\subsection{Discrete Green functions}\quad Our first result is for the circulant matrices $A_1$ and $A_2$.
\begin{lemma}\label{tn:algebra}
The matrices $A_1$ and $A_2$ commute with each other, and satisfy
\begin{equation}\label{tn:algebraeq}
 A_1^T=-A_1,\quad A_2^T=A_2,\quad
 A_1\mathbf1=A_2\mathbf1=\mathbf 0\quad\text{and}\quad
 A_1^2=A_2(A_2+4I).
\end{equation}
Moreover,
\begin{equation*}\label{tn:dirichletform}
 \mathbf a^TA_2\mathbf a=-\sum_{j=0}^{K-1}(a_{j+1}-a_j)^2
\end{equation*}
and the inverse $A_2^{-1}$ is well defined on $\mathbf1^\perp$ and satisfies
\begin{equation}\label{tn:eliminationidentity}
 A_1A_2^{-1}A_1\mathbf a=(A_2+4I)\mathbf a
 \quad\text{for}\quad\bar a=0.
\end{equation}
\end{lemma}
\begin{proof}
Let $S\mathbf a=(a_{1}, a_2, \cdots, a_{K-1}, a_0)$ for a cyclic sequence $\mathbf a=(a_0, a_1, \cdots, a_{K-1})$. Then $S^{-1}=S^T$ with
\[
 A_1=S-S^{-1}\quad\text{and}\quad A_2=S+S^{-1}-2I.
\]
Thus, $A_1$ and $A_2$ are polynomials in $S$ and $S^{-1}$.  It follows that $ A_1^T=-A_1$, $A_2^T=A_2$ and $A_1\mathbf1=A_2\mathbf1=\mathbf 0$.  Moreover, by direct calculations,
\begin{eqnarray*}
A_2(A_2+4I)&=&(S+S^{-1}-2I)(S+S^{-1}+2I)\\
 &=&S^2+2I+S^{-2}-4I\\
 &=&A_1^2
\end{eqnarray*}
and
\begin{align*}
 \sum_{j=0}^{K-1}a_j(A_2a)_j=2\sum_{j=0}^{K-1}a_ja_{j+1}-2\sum_{j=0}^{K-1}a_j^2
 =-\sum_{j=0}^{K-1}(a_{j+1}-a_j)^2.
\end{align*}
Consequently, $\ker A_2=\operatorname{span}\{\mathbf1\}$ and
$\operatorname{ran}A_2=\mathbf1^\perp$.  Thus, for a cyclic sequence $\mathbf a=(a_0, a_1, \cdots, a_{K-1})$ satisfying $\bar a=0$, commutation and
\eqref{tn:algebraeq} give
$$
A_1A_2^{-1}A_1\mathbf a=A_1^2A_2^{-1}\mathbf a=(A_2+4I)\mathbf a,
$$
which completes the proof.
\end{proof}

For $0\le l<K$, we introduce the discrete Fourier transform and its inverse for cyclic sequences
\begin{equation}\label{tn:fourier}
\widehat a_l=\frac1K\sum_{j=0}^{K-1}a_j e^{-\mathrm i ljh}\quad\text{and}\quad
 a_j=\sum_{l=0}^{K-1}\widehat a_l e^{\mathrm i ljh},
\end{equation}
respectively.
The geometric series identity gives
\begin{eqnarray}\label{identity}
 \frac1K\sum_{j=0}^{K-1}e^{\mathrm i(l-r)jh}=\delta_{lr}\quad\text{and}\quad
 \|\mathbf a\|_2^2=K\sum_{l=0}^{K-1}|\widehat a_l|^2.
\end{eqnarray}
Moreover, substitution into the two differences $A_1$ and $A_2$ yields
\begin{equation}\label{tn:symbols}
 \widehat{(A_1{\bf a})_l}=2\mathrm i\sin(lh)\widehat a_l\quad\text{and}\quad
 \widehat{(A_2{\bf a})_l}=-4s^2t_{l,K}\widehat a_l,
\end{equation}
where $t_{l,K}=\frac{\sin^2(l\pi/K)}{s^2}$ with $s=\sin\frac{\pi}{K}$.
Clearly, we have $t_{l,K}\ge1$ for $l\ne0$.  If $r_{l,K}=\min\{l,K-l\}$ then concavity of the function $\sin(x)$ also gives
\begin{equation}\label{tn:frequencies}
 t_{l,K}\ge\frac{4r_{l,K}^2}{\pi^2}\quad\text{and}\quad
 \sum_{l=1}^{K-1}\frac1{t_{l,K}}\le
 \frac{\pi^2}{2}\sum_{r=1}^\infty\frac1{r^2}<\infty.
\end{equation}

\vskip0.12in

We recall that in applications, we have
\begin{equation*}\label{tn:nuexact}
 \nu_{*,K}=
 \frac{(m+1)a_0mR_K^{-m-1}+(n+1)b_0nR_K^{-n-1}}
      {a_0mR_K^{-m-1}+b_0nR_K^{-n-1}}
\end{equation*}
by \eqref{eq:D-nueff} and Lemma~\ref{pf:radius}, which implies $\nu_{*,K}\to m_0+1>3$ and $\lambda_{*,K}\to m_0-1>1$
as $K\to\infty$, where
\begin{equation*}\label{tn:lambda}
\lambda_{*,K}=(\nu_{*,K}-2)c^2-s^2+\frac{c^2+\nu_{*,K}s^2}{D_K}
\end{equation*}
with $c=\cos\frac{\pi}{K}$ and $D_K=d_K+O(1)\to+\infty$ as $K\to\infty$.
Thus, we may assume that
\begin{equation}\label{tn:parameterwindow}
2+\epsilon_0\le\nu_{*,K}\le M\quad\text{and}\quad D_K\ge D_0
\quad\text{for}\quad K\ge K_0,
\end{equation}
and
\begin{equation}\label{tn:lambdabounds}
0<\lambda_-\le\lambda_{*,K}\le\lambda_+<\infty\quad\text{for}\quad
D_K-\nu_{*,K}\asymp D_K,
\end{equation}
by increasing $K_0$ if necessary in terms of $\epsilon_0$ and $M$, where $\lambda_-$ and $\lambda_+$ are positive constants depending only on $\epsilon_0$ and $M$.
\begin{proposition}
The matrix $T_\parallel$, given by \eqref{eq:matrixblocks}, has $K-1$ negative eigenvalues $\lambda_{l,K}^-$, one zero eigenvalue and $K$ positive eigenvalues $\lambda_{l,K}^+$ with $\lambda_{0,K}^+=D_K-\nu_{*,K}$,
\begin{equation}\label{tn:eigenbounds}
-\lambda_{l,K}^-\asymp t_{l,K}\quad\text{and}\quad
\lambda_{l,K}^+\asymp D_K(1+t_{l,K})
\end{equation}
for $0\leq l\leq K-1$ as $K\to\infty$.  In particular,
\begin{equation}\label{tn:planarkernel}
 \ker T_\parallel=\operatorname{span}\{(\mathbf0,\mathbf1)\}
\end{equation}
and if $(\mathbf f,\mathbf g)$ solves the equation
\begin{eqnarray}\label{eqn0001}
\left\{\aligned
&T_\parallel\binom{\mathbf f}{\mathbf g}=(\mathbf a,\mathbf b),\\
&\bar g=\overline{b}=0,
\endaligned\right.
\end{eqnarray}
then we have
\begin{equation}\label{tn:l2estimate}
 \|\mathbf f\|_2+\|\mathbf g\|_2+
 \|\dot{\mathbf f}\|_2+\|\dot{\mathbf g}\|_2+
 \|\ddot{\mathbf f}\|_2+\|\ddot{\mathbf g}\|_2
 \lesssim(\|\mathbf a\|_2+\|\mathbf b\|_2).
\end{equation}
\end{proposition}
\begin{proof}
The computation for the eigenvalues and the kernel of $T_\parallel$ is contained in \cite[Lemma~8.3]{DelPinoWeiYao2015}.  However, we use the Fourier transform to recompute it and improve the estimates in \cite[Lemma~8.3]{DelPinoWeiYao2015} to \eqref{tn:l2estimate}.  More precisely, by \eqref{tn:symbols}, the $l$-th matrix of the Fourier transform of $T_\parallel$ is given by
\begin{equation*}\label{tn:modeblock}
T_\parallel(l)=\begin{pmatrix}
 D_K-\nu_{*,K}-(c^2+D_Ks^2)t_{l,K}&2\mathrm i b_{e,K}\sin(lh)\\
 -2\mathrm i b_{e,K}\sin(lh)&(D_Kc^2+s^2)t_{l,K}
 \end{pmatrix},
\end{equation*}
where $b_{e,K}$ is given by \eqref{eq:matrixcoefficients} with $D_K$ given by Lemma~\ref{pf:radius}.
Thus
\begin{equation}\label{tn:determinant}
\det T_\parallel(l)=-D_Kt_{l,K}(t_{l,K}+\lambda_{*,K})
\end{equation}
for all $0\leq l\leq K-1$,
which implies the complete information on the eigenvalues and the kernel of $T_\parallel$.
The estimate~\eqref{tn:eigenbounds} then follows from direct calculations under the conditions~\eqref{tn:parameterwindow} and \eqref{tn:lambdabounds}.  By \eqref{tn:determinant}, we also have $\|T_\parallel(l)^{-1}\|_2\lesssim t_{l,K}^{-1}$, where $\|T_\parallel(l)^{-1}\|_2$ is the Euclidean norm of the matrix $T_\parallel(l)^{-1}$.  Thus, for a solution of \eqref{eqn0001},
\begin{eqnarray*}
\sum_{l=0}^{K-1}\binom{\widehat{a}_l}{\widehat{b}_l}e^{\mathrm i ljh}=\sum_{l=0}^{K-1}T_\parallel(l)\binom{\widehat{f}_l}{\widehat{g}_l}e^{\mathrm i ljh}.
\end{eqnarray*}
Since
\[
\widehat{\dot f}_l
=\frac{e^{\mathrm i lh}-1}{h}\widehat f_l,
\quad\widehat{\ddot f}_l=-\frac{4s^2t_{l,K}}{h^2}\widehat f_l
\]
and
\[
\widehat{\dot g}_l
=\frac{e^{\mathrm i lh}-1}{h}\widehat g_l,
\quad
\widehat{\ddot g}_l=-\frac{4s^2t_{l,K}}{h^2}\widehat g_l,
\]
by $\|T_\parallel(l)^{-1}\|_2\lesssim t_{l,K}^{-1}$, we have
\begin{align}
|\widehat{\dot f}_l|^2+|\widehat{\dot g}_l|^2
\le t_{l,K}\bigl(|\widehat f_l|^2+|\widehat g_l|^2\bigr)\lesssim\frac {1}{t_{l,K}}
\bigl(|\widehat{a}_l|^2+|\widehat{b}_l|^2\bigr)\lesssim\bigl(|\widehat{a}_l|^2+|\widehat{b}_l|^2\bigr)
\label{eq:first-bound}
\end{align}
and
\begin{align}
|\widehat{\ddot f}_l|^2+|\widehat{\ddot g}_l|^2
\le t_{l,K}^2\bigl(|\widehat f_l|^2+|\widehat g_l|^2\bigr)\lesssim\bigl(|\widehat{a}_l|^2+|\widehat{b}_l|^2\bigr).
\label{eq:second-bound}
\end{align}
Thus, by \eqref{identity}, \eqref{eq:first-bound} and \eqref{eq:second-bound},
\begin{eqnarray*}
&&\|\mathbf f\|_2^2+\|\mathbf g\|_2^2
 +\|\dot{\mathbf f}\|_2^2+\|\dot{\mathbf g}\|_2^2
 +\|\ddot{\mathbf f}\|_2^2+\|\ddot{\mathbf g}\|_2^2\notag\\
&=&K\sum_{l=0}^{K-1}
 \bigl(|\widehat f_l|^2+|\widehat g_l|^2
 +|\widehat{\dot f}_l|^2+|\widehat{\dot g}_l|^2
 +|\widehat{\ddot f}_l|^2+|\widehat{\ddot g}_l|^2\bigr)\notag\\
&\lesssim& K\sum_{l=0}^{K-1}
 \bigl(|\widehat{a}_l|^2+|\widehat{b}_l|^2\bigr)
\lesssim\bigl(\|\mathbf a\|_2^2+\|\mathbf b\|_2^2\bigr).
\label{eq:sum-of-squares}
\end{eqnarray*}
It completes the proof.
\end{proof}

To pass to the norm used in the reduction, we introduce
\begin{eqnarray*}
\mathbf G_K=(G_K(0), \cdots, G_K(j), \cdots, G_K(K-1)),
\end{eqnarray*}
where
\begin{equation}\label{tn:poissongreen}
 G_K(j)=\frac{K^2-1}{12K}-\frac{j(K-j)}{2K}\quad \text{for}\quad 0\le j\leq K-1
\end{equation}
and extend it periodically to $\mathbb{Z}$.  This is the classical discrete Green function of the cycle, see
\cite{ChungYau2000} for the general theory and \cite[Theorem~1]{Ellis2003} for the explicit
formula (the factor $2$ reflects the Laplacian normalization).  More precisely, $\mathbf G_K$ is the discrete Green function of $A_2$ satisfying
\begin{equation}\label{tn:greenidentity}
\left\{\aligned
&-A_2\mathbf G_K=\delta_0-K^{-1}\mathbf1,\\
&\sum_{j=0}^{K-1}G_K(j)=0,
\endaligned\right.
\end{equation}
where $\delta_0$ is the Kronecker delta given by
\begin{eqnarray*}
\delta_0(j)=\left\{\aligned
1,\quad j\equiv0\ \ (\text{mod }K),\\
0,\quad j\not\equiv0\ \ (\text{mod }K).
\endaligned\right.
\end{eqnarray*}
\begin{lemma}\label{tn:greenpoisson}
The unique solution of the following equation
\begin{equation}\label{tn:greenidentity-1}
\left\{\aligned
&-A_2\mathbf f=\mathbf a,\\
&\sum_{j=0}^{K-1}f(j)=0,
\endaligned\right.
\end{equation}
for $\mathbf a$ satisfying $\bar a=0$, is given by
\begin{equation}\label{tn:poissonsolution}
 f_j=\sum_{i=0}^{K-1}G_K(j-i)a_i.
\end{equation}
It satisfies
\begin{equation*}\label{tn:poissonestimate}
 \|\mathbf f\|_*\lesssim K^2\|\mathbf a\|_\infty\quad\text{and}\quad
 \|A_1A_2^{-1}\mathbf a\|_\infty\lesssim K\|\mathbf a\|_\infty.
\end{equation*}
\end{lemma}
\begin{proof}
Since $\ker A_2=\operatorname{span}\{\mathbf1\}$ and
$\operatorname{ran}A_2=\mathbf1^\perp$ by Lemma~\ref{tn:algebra} and $\mathbf G_K$ is the discrete Green function of $A_2$ satisfying \eqref{tn:greenidentity}, by the discrete convolution formula,
\begin{eqnarray*}
\mathbf f=(-A_2)^{-1}\mathbf a=G_K*\mathbf a
\end{eqnarray*}
with $f_j$ given by \eqref{tn:poissonsolution} for the solution of \eqref{tn:greenidentity-1} and $\mathbf a$ satisfying $\bar a=0$.  The explicit polynomial of \eqref{tn:poissongreen} gives
\[
 \max_{0\leq j\leq K-1}|G_K(j)|\lesssim K\quad\text{and}\quad
 \max_{0 \leq j \leq K-1}|G_K(j+1)-G_K(j)|\lesssim1.
\]
Moreover, $Ks\le\pi$ and $s\asymp h$.  Thus, $\|\mathbf f\|_*\lesssim K^2\|\mathbf a\|_\infty$ by
\begin{eqnarray*}
\left\{\aligned
&\|\mathbf f\|_\infty
\lesssim K^2\|\mathbf a\|_\infty,\,\,\,\,
\|\dot{\mathbf f}\|_\infty\lesssim \frac{K}{h}\|\mathbf a\|_\infty\lesssim K^2\|\mathbf a\|_\infty,\\
 &\|\ddot{\mathbf f}\|_\infty=\frac{1}{h^2}\|\mathbf a\|_\infty\lesssim\frac{1}{h^2}\|\mathbf a\|_\infty\lesssim K^2\|\mathbf a\|_\infty
 \endaligned\right.
\end{eqnarray*}
and we also have
\[
\|A_1A_2^{-1}\mathbf a\|_\infty\lesssim\max_{0 \leq j \leq K-1}\sum_{i=0}^{K-1}|G_K(j+1-i)-G_K(j-1-i)|\left\|a_i\right\|_\infty\lesssim K\|\mathbf a\|_\infty,
\]
which completes the proof.
\end{proof}

To analyze the matrix~$T_N$, we also need the following.
\begin{lemma}\label{tn:resolvent}
For $\lambda\in[\lambda_-,\lambda_+]\subset(0,\infty)$, where $\lambda_\pm$ are given by \eqref{tn:lambdabounds}, the positive Green sequence of $-\frac{A_2}{4s^2}+\lambda I$ is given by
$$
\mathbf G_{K,\lambda}=(G_{K,\lambda}(0), G_{K,\lambda}(1), \cdots, G_{K,\lambda}(K-1))
$$
with
\begin{equation}\label{tn:resolventgreen}
 G_{K,\lambda}(j)=\frac{2s^2}{\sinh\kappa}
 \frac{e^{-\kappa j}+e^{-\kappa(K-j)}}{1-e^{-\kappa K}}
\end{equation}
for $0\le j<K$ and extend it periodically to $\mathbb{Z}$, where $\kappa=\text{arccosh}(1+2s^2\lambda)$ and $\sum_{j=0}^{K-1}G_{K,\lambda}(j)=\lambda^{-1}$.  Moreover, the unique solution of the equation
\begin{equation*}\label{tn:positiveresolvent}
 \left(-\frac{A_2}{4s^2}+\lambda I\right)\mathbf f=\mathbf a
\end{equation*}
is given by $f_j=\sum_iG_{K,\lambda}(j-i)a_i$ and satisfies
\begin{equation*}\label{tn:resolventbound}
 \|\mathbf f\|_*\lesssim\|\mathbf a\|_\infty.
\end{equation*}
\end{lemma}
\begin{proof}
The closed form \eqref{tn:resolventgreen} is the cyclic one-dimensional lattice Green
function, see \cite[Eq.~II.19]{MontrollWeiss1965}, \cite{Feller1971},
\cite{KatsuraEtAl1971} and \cite[Eqs.~(42)--(45)]{Dowker2012} in the present periodized form
with $\sinh^2\gamma=s^2\lambda$, $\kappa=2\gamma$, or \cite[Theorem~8 and Corollary~5]{Ellis2003}
with $\alpha=2s^2\lambda$, $r=e^{\kappa}$.  Positivity follows from the random-walk
(Neumann series) representation, see also \cite{AldousFill2002}.  More precisely, it is easy to check that
\begin{equation*}\label{eq:step1}
A_2g(j)
=(e^{-\kappa}+e^{\kappa}-2)e^{-\kappa j}\quad\text{for}\quad j\geq1,
\end{equation*}
where $g(j)=e^{-\kappa|j|}$.  Thus, by taking $\kappa=arccosh(1+2s^2\lambda)$, we have
\[
-\frac{A_2g(j)}{4s^2}+\lambda g(j)=\frac{\sinh\kappa}{2s^2}\delta_0.
\]
It follows from \eqref{tn:resolventgreen} that the $K$-periodic extension of
$$
\mathbf G_{K,\lambda}=(G_{K,\lambda}(0), G_{K,\lambda}(1), \cdots, G_{K,\lambda}(K-1))
$$
is the positive Green sequence of $-\frac{A_2}{4s^2}+\lambda I$ with $\lambda\sum_{j=0}^{K-1}G_{K,\lambda}(j)=1$.
By the discrete convolution formula,
\begin{eqnarray*}
\mathbf f=\left(-\frac{A_2}{4s^2}+\lambda I\right)^{-1}\mathbf a=\mathbf G_{K,\lambda}*\mathbf a
\end{eqnarray*}
with $f_j=\sum_iG_{K,\lambda}(j-i)a_i$ and
\begin{eqnarray*}
\|\mathbf f\|_\infty\lesssim\sum_{j=0}^{K-1}G_{K,\lambda}(j)\|\mathbf a\|_\infty\lesssim\|\mathbf a\|_\infty.
\end{eqnarray*}
Moreover, by $\lambda\in[\lambda_-,\lambda_+]\subset(0,\infty)$, $\kappa\asymp s\asymp h$ and the geometric series identity,
\begin{eqnarray*}
\|\dot{\mathbf f}\|_\infty\lesssim\left(\sum_{j=0}^{K-1}|\dot{G}_{K,\lambda}(j)|\right)\|\mathbf a\|_\infty\lesssim\frac{s^2}{h\kappa}\|\mathbf a\|_\infty\lesssim\|\mathbf a\|_\infty
\end{eqnarray*}
and by $(-\frac{A_2}{4s^2}+\lambda I)\mathbf G_{K,\lambda}=\delta_0$ with $\lambda\sum_{j=0}^{K-1}G_{K,\lambda}(j)=1$,
\begin{eqnarray*}
\|\ddot{\mathbf f}\|_\infty\lesssim\sum_{j=0}^{K-1}|\ddot{G}_{K,\lambda}(j)|\|\mathbf a\|_\infty\lesssim\sum_{j=0}^{K-1}\frac{|A_2G_{K,\lambda}(j)|}{h^2}\|\mathbf a\|_\infty\lesssim\|\mathbf a\|_\infty,
\end{eqnarray*}
which completes the proof.
\end{proof}

\subsection{Inverse of $T_N$}\quad We are now in the position to analyze the inverse of $T_N$ for $K>0$ sufficiently large.  We begin with the analysis of the inverse of $T_\parallel$.
\begin{proposition}\label{tn:planarinverse}
Under the condition~\eqref{tn:parameterwindow}, the system
\begin{equation}\label{tn:system}
\begin{cases}
 a_{e,K}A_2\mathbf f+(D_K-\nu_{*,K})\mathbf f+b_{e,K}A_1\mathbf g=\mathbf a,\\
 -b_{e,K}A_1\mathbf f+c_{e,K}A_2\mathbf g=\mathbf b,
\end{cases}
\end{equation}
is solvable if and only if $\overline{b}=0$, where $a_{e,K}$, $b_{e,K}$ and $c_{e,K}$ are given by \eqref{eq:matrixcoefficients} with $D_K$ given by Lemma~\ref{pf:radius}.  Moreover, its solution with $\bar g=0$
is unique and satisfies
\begin{equation}\label{tn:planarstar}
 \|\mathbf f\|_*+\|\mathbf g\|_*
 \lesssim (\|\mathbf a\|_\infty+\|\mathbf b\|_\infty).
\end{equation}
\end{proposition}
\begin{proof}
This is just the discrete
counterpart of
\cite[Lemma~8.4]{DelPinoWeiYao2015}.  By \eqref{tn:differences}, summing both equations in \eqref{tn:system}, respectively, gives
\begin{equation*}\label{tn:means}
 (D_K-\nu_{*,K})\bar f=\overline{a}\quad\text{and}\quad \overline{b}=0.
\end{equation*}
Thus $\bar f=\frac{\overline{a}}{D_K-\nu_{*,K}}$ is fixed, and by Lemma~\ref{tn:algebra}, \eqref{tn:system} can be solved only if $\bar b=0$.  Moreover, under the orthogonal condition~$\bar g=0$ and by Lemma~\ref{tn:algebra}, the second equation of \eqref{tn:system} then gives the exact recovery formula
\begin{equation}\label{tn:recoverg}
 \mathbf g=c_{e,K}^{-1}A_2^{-1}\mathbf b
      +\frac{b_{e,K}}{c_{e,K}}A_2^{-1}A_1(\mathbf f-\bar f\mathbf1)
\end{equation}
for $\overline{b}=0$.
Applying $A_1$ to \eqref{tn:recoverg} and substituting in the first equation of \eqref{tn:system}, together with \eqref{tn:eliminationidentity}, gives
\begin{align*}
\left[\left(a_{e,K}+\frac{b_{e,K}^2}{c_{e,K}}\right)A_2
       +\left(D_K-\nu_{*,K}+\frac{4b_{e,K}^2}{c_{e,K}}\right)I\right]
       (\mathbf f-\bar f\mathbf1)=\mathbf a-\overline{a}\mathbf1
                         -\frac{b_{e,K}}{c_{e,K}}A_1A_2^{-1}\mathbf b,
\end{align*}
which, together with \eqref{eq:matrixcoefficients}, implies that
\begin{align}\label{tn:schurequation}
\left(-\frac{A_2}{4s^2}+\lambda_{*,K}I\right)
                (\mathbf f-\bar f\mathbf1)=-\frac{D_Kc^2+s^2}{D_K}
              (\mathbf a-\overline{a}\mathbf1)
     -\frac{(D_K-1)c}{D_K}sA_1A_2^{-1}\mathbf b.
\end{align}
Since every term on the right side of \eqref{tn:schurequation} has zero mean by $\overline{b}=0$, by Lemma~\ref{tn:resolvent}, \eqref{tn:schurequation} is unique solvable and the unique solution $\mathbf f$ satisfies
\begin{equation}\label{tn:fstar}
 \|\mathbf f\|_*\lesssim(\|\mathbf a\|_\infty+\|\mathbf b\|_\infty).
\end{equation}
It follows from \eqref{tn:recoverg} that \eqref{tn:system} can be solved if and only if $\bar b=0$, and the solution is unique under the orthogonal condition~$\bar g=0$.  Moreover, since by \eqref{eq:matrixcoefficients}, $|c_{e,K}|\asymp D_K/s^2\asymp K^2\log K$ and $|b_{e,K}/c_{e,K}|\lesssim s\asymp K^{-1}$ as $K\to+\infty$.  It follows from Lemma~\ref{tn:greenpoisson} that
\begin{equation}\label{tn:gstar}
 \|\mathbf g\|_*\lesssim(\|\mathbf a\|_\infty+\|\mathbf b\|_\infty).
\end{equation}
\eqref{tn:planarstar} follows from \eqref{tn:fstar} and \eqref{tn:gstar}, which completes the proof.
\end{proof}

We next move to the analysis of the inverse of $T_\perp=I+\frac{A_2}{4s^2}$.
\begin{proposition}\label{tn:normalinverse}
$T_\perp$, i. e., the normal block in $T_N$, has one positive eigenvalue, two zero eigenvalues and
$K-3$ negative eigenvalues with the kernel
\begin{equation}\label{tn:normalkernel}
 \ker T_\perp=\operatorname{span}
       \{(\cos\theta_j)_j,(\sin\theta_j)_j\}.
\end{equation}
Moreover, if $\mathbf a$ is orthogonal to $ \ker T_\perp$, then the equation
$T_\perp\mathbf z=\mathbf a$ has a unique solution which is also orthogonal to $ \ker T_\perp$ with
\begin{equation}\label{tn:normalstar}
 \|\mathbf z\|_*\lesssim \|\mathbf a\|_\infty.
\end{equation}
\end{proposition}
\begin{proof}
By~\eqref{tn:symbols}, after applying the Fourier transform on $T_\perp$, the eigenvalue on frequency $l$ is $1-t_{l,K}$, where $t_{l,K}$ is given by \eqref{tn:fourier}.  Thus, $t_{0,K}=0$, $t_{1,K}=t_{K-1,K}=1$ and $t_{l,K}>1$ for $2\leq l\leq K-2$.  It follows that $T_\perp$ has one positive eigenvalue, two zero eigenvalues and
$K-3$ negative eigenvalues with the kernel given by \eqref{tn:normalkernel}.  By \eqref{tn:fourier}, the exact solution of $T_\perp\mathbf z=\mathbf a$ is
\[
 z_j=\sum_{l=0\atop
 l\not=1,K-1}^{K-1}\frac{\widehat{a}_l}{1-t_{l,K}}e^{\mathrm iljh}
\]
for $\mathbf a\in(\ker T_\perp)^{\perp}$,
where $\widehat{a}_l$ is the Fourier coefficients of $\bf a$ given by \eqref{tn:fourier}.  Thus, $\bf z$ is also orthogonal to $ \ker T_\perp$.  Let $r_{l,K}=\min(l,K-l)\ge2$, then by \eqref{tn:frequencies},
\[
\frac{1}{|1-t_{l,K}|}\lesssim\frac{1}{r_{l,K}^2}\quad\text{for}\quad r_{l,K}\ge3,
\]
which, together with \eqref{identity}, implies that $\|\mathbf z\|_\infty\lesssim\|\mathbf a\|_\infty$.  By $T_\perp=I+\frac{A_2}{4s^2}$,
\[
\|\ddot{\mathbf z}\|_\infty\lesssim\left\|\frac{A_2\mathbf z}{h^2}\right\|_\infty\lesssim\left\|\frac{4s^2}{h^2}(\mathbf a-\mathbf z)\right\|_{\infty}\lesssim\|\mathbf a\|_\infty.
\]
By the discrete interpolation inequality, we also have
\[
\|\dot{\mathbf z}\|_\infty\lesssim(\|\ddot{\mathbf z}\|_\infty\|\mathbf z\|_\infty)^{\frac12}\lesssim\|\mathbf a\|_\infty,
\]
which implies that \eqref{tn:normalstar} holds true.
\end{proof}

By \eqref{tn:planarkernel} and \eqref{tn:normalkernel}, the kernel of $T_N$ is spanned by
$\mathbf k_1=(\mathbf0,\mathbf1)$ for $N=2$ and
\begin{eqnarray}\label{ker}
 \mathbf k_1=(\mathbf0,\mathbf1,\mathbf0),\quad
 \mathbf k_2=(\mathbf0,\mathbf0,(\sin\theta_j)_j)\quad\text{and}\quad
 \mathbf k_3=(\mathbf0,\mathbf0,(-\cos\theta_j)_j)
\end{eqnarray}
for $N=3$.  Moreover, it is easy to see that
\begin{eqnarray}\label{orth-ker}
 \mathbf k_j\cdot\mathbf k_l=0\ (j\ne l)\quad\text{and}\quad
 |\mathbf k_1|^2=2|\mathbf k_2|^2=2|\mathbf k_3|^2=K.
\end{eqnarray}
Let $\mathcal E_K=(\ker T_N)^\perp$ and $\Pi_{\mathcal E_K}: \R^{NK}\to (\ker T_N)^\perp$ be
\begin{equation*}\label{tn:projection}
 \Pi_{\mathcal E_K}\mathbf a=\mathbf a-\sum_{j=1}^{N(N-1)/2}
\frac{\mathbf a\cdot\mathbf k_j}{|\mathbf k_j|^2}\mathbf k_j.
\end{equation*}
\begin{proposition}\label{tn:TNestimate}
For every $\mathbf a\in\mathcal E_K$, there is a unique
$\overline{\mathbf q}\in\mathcal E_K$ such that
$T_N\overline{\mathbf q}=\mathbf a$ with
\begin{equation}\label{tn:TNstar}
 \|\overline{\mathbf q}\|_*\lesssim\|\mathbf a\|_\infty.
\end{equation}
\end{proposition}
\begin{proof}
By \eqref{s:eq:T3} and Propositions~\ref{tn:planarinverse} and \ref{tn:normalinverse},
$T_N\overline{\mathbf q}=\mathbf a$ is unique solvable for $\mathbf a\in\mathcal E_K$ with the estimate~\eqref{tn:TNstar}.
\end{proof}

\end{document}